\documentclass[a4paper,12pt]{amsart}

\usepackage{graphicx}
\usepackage{subfigure}
\usepackage[usenames, dvipsnames]{color}

\newcommand{\RR}{{\bf R}}

\newcommand{\ZZ}{{\bf Z}}

\newcommand{\xy}{(x,y)}

\newcommand{\seta}{\longrightarrow}
\newcommand{\dpt}{\displaystyle}
\newcommand{\va}{\varphi}
\newcommand{\ba}{{\bf a}}
\newcommand{\FF}{{\mathcal F}}
\newcommand{\GG}{{\mathcal G}}

\DeclareMathOperator{\sign}{sign}

\newtheorem{theorem}{Theorem}[section]
\newtheorem{lemma}[theorem]{Lemma}
\newtheorem{proposition}[theorem]{Proposition}
\newtheorem{corollary}[theorem]{Corollary}

\theoremstyle{definition}
\newtheorem{Ex}[theorem]{Example}

\begin{document}
\title[Qualitative planar dynamics with star nodes: \today]{Qualitative planar dynamics with star nodes and homogeneous nonlinearities: beyond Hilbert's $16^{th}$ problem\\ \today}

\author[B.Alarc\'on]{Bego\~na Alarc\'on}
\address{B.Alarc\'on\\Departamento de Matem\'atica Aplicada, Instituto de Matem\'atica e Estat\'\i stica, Universidade Federal Fluminense \\
Rua Professor Marcos Waldemar de Freitas Reis, S/N, Bloco H \\
Campus do Gragoat\'a,
CEP 24.210 -- 201, S\~ao Domingos --- Niter\'oi, RJ, Brasil}
\email{balarcon@id.uff.br}
\thanks{The first author was partially supported by the Spanish Research Project PID2020--113052GB--I00}
\author[S.B.S.D. Castro]{Sofia B.S.D. Castro}
\address{S.B.S.D. Castro\\Centro de Matem\'atica da Universidade do Porto\\ Rua do Campo Alegre 687, 4169-007 Porto, Portugal
and
Faculdade de Economia do Porto \\ Rua Dr. Roberto Frias, 4200-464 Porto, Portugal}
\email{sdcastro@fep.up.pt}
\thanks{The last two authors were partially supported by CMUP, member of LASI, which is financed by national  funds through  FCT --- Funda\c c\~ao para a Ci\^encia e a Tecnologia, I.P. (Portugal)  under the project with reference  UID/00144/2025.}
\author[I.S. Labouriau]{Isabel S. Labouriau}
\address{I.S. Labouriau\\Centro de Matem\'atica da Universidade do Porto\\ Rua do Campo Alegre 687, 4169-007 Porto, Portugal}
\email{islabour@fc.up.pt}

\begin{abstract}
This is a  full study of the dynamics of polynomial planar vector fields whose linear part is a multiple of the identity and whose nonlinear part is a homogeneous polynomial of arbitrary degree $n>1$.
It extends previous work by other authors that was mainly concerned with the existence and number of limit cycles.
The general results are also applied to two classes of examples where the nonlinearities have degrees 2 and 3,
for which we provide a  set of phase portraits.
\end{abstract}

\maketitle



%
%
%

\section{Introduction}

Global planar dynamics of polynomial vector fields has been of interest for many years. Part of this interest arises from its connection to Hilbert's $16^{th}$ problem on the number of limit cycles for the dynamics.
Because of Hilbert's $16^{th}$ problem, substantial effort has been devoted to establishing a bound for the number of limit cycles. For some contributions in this direction when the vector field has homogeneous nonlinearities see the work of Boukoucha \cite{B2017}, Bendjeddou {\em et al.} \cite{BLS2013}, Huang {\em et al.} \cite{HLL}, Gasull {\em et al.} \cite{GYZ}, Llibre {\em et al.}  \cite{LYZ} or Carbonell and Llibre \cite{CL1989}, and more recently Garc\'{\i}a-Salda\~na {\em et al.} \cite{Gar-SalGasGia}. This question has also been approached using bifurcations by, for instance, Benterki and Llibre \cite{BL} or \cite{GYZ}.
 Vector fields of degree 2 are reviewed in the book by Art\'es {\em et al.} \cite{Artes}.
Quasi-homogeneous nonlinearities have recently been addressed by Llibre {\em et al.} \cite{Llibre2025}.
The special case where the   homogeneous nonlinearity is contracting is treated in Alarc\'on {\em et al.} \cite{Alarcon}.
Problems with symmetry appear in \'Alvarez {\em et al.} \cite{AGP} and Labouriau and Murza \cite{LM}.\footnote{Our references do not pretend to be comprehensive. The reader can find further interesting work by looking at the references within those we provide.}

We are, of course, also concerned in establishing an upper bound for the number of limit cycles. However, when no limit cycle exists, we take a different route and address the  question of the existence of polycycles\footnote{Including  heteroclinic cycles, that are  particular instances of  polycycles.} and the number of equilibria in them. 
Our main results are stated and proved in Sections~\ref{sec:equilibria} -- \ref{sec:regions}.

We focus on polynomial vector fields with homogeneous nonlinearities, as many before us, but use the existence of invariant 
 lines through the origin to provide information on the global dynamics. We consider vector fields of the form 
 \begin{equation} \label{eqR2contracting}
\left\{\begin{array}{lcl}
\dot{x}&=& \lambda x +Q_1(x,y)  \\
\mbox{} & \mbox{} & \\
\dot{y} &=& \lambda y +Q_2(x,y) ,
\end{array}  \right. 
\end{equation}
where $\lambda \neq 0$ and $Q_i$, $i=1,2$  are homogeneous non-zero polynomials
 of the same degree  $n>1$.
 We define $Q=(Q_1, Q_2)$ and say it is a homogeneous polynomial of degree $n$.
  The origin of such a system is an  {\em unstable star node} (a node with equal and positive eigenvalues) if $\lambda >0$ and a  {\em stable star node} (a node with equal and negative eigenvalues) if $\lambda <0$.

Using both the dynamics of the vector field at infinity, through its Poincar\'e compactification, as well as polar coordinates, we describe completely the existence of equilibria at infinity. These occur as equilibria on the boundary of the Poincar\'e disk for the dynamics of the compactification. We distinguish them from finite equilibria, occurring in the interior of the Poincar\'e disk, by calling them {\em infinite equilibria}. The existence of infinitely many infinite equilibria determines that either there are also infinitely many finite equilibria or the origin is the only finite equilibrium. 

The case of finitely many infinite equilibria allows for a complete description of the planar dynamics. Each equilibrium at infinity defines an invariant radius with the origin. All finite equilibria are located on such radii. We are able to provide an upper bound for the number of finite equilibria which improves on Bezout's Theorem for polynomials of degree strictly greater than 3. We also describe the stability of all equilibria.

The above is a preliminary step for our main contribution to the description of the global planar dynamics with star nodes. This relies on the construction of invariant sectors from the invariant radii. The number of sectors necessary for the full description of the global dynamics depends on the number of infinite equilibria (whose upper bound we have established previously). Polycycles are present only when one type of sector repeats to cover the plane. A limit cycle exists only when there are no infinite equilibria. 
It can arise as a generic perturbation of a heteroclinic cycle that collapses through a saddle-node bifurcation.
Together with ours, results previously established by Bendjeddou {\em et al.} \cite{BLS2013}, Coll {\em et al.} \cite{CGP1997}, and Gasull {\em et al.} \cite{GYZ} completely describe the case when the origin is a global attractor or repellor.
We describe the possible phase portraits up to {\em topological equivalence}, where two systems  \eqref{eqR2contracting} are equivalent if there is a homeomorphism of the Poincar\'e disk mapping trajectories of one system into trajectories of the other and preserving time orientation.

We illustrate our results by  studying the cases when the nonlinear polynomials are of degree 2 and 3. To this purpose, we build also on previous results of Cima and Llibre \cite{AnnaL1990} and Date \cite{Date}.  
In both cases we provide a list of admissible phase portraits for the global dynamics.
In the case of degree 2 this can be complemented with the classification  in the book by Art\'es {\em et al.} \cite{Artes}.
This article is thus a preliminary step towards the study of the probability of occurrence of a given phase portrait along the lines of Cima {\em et al.} \cite{CimGasMan}.
For higher degree the interested reader may apply the same procedure to the classification of Collins \cite{Collins}.

\bigbreak

This article is organised as follows: in the next section we provide some background and establish our notation. In Section~\ref{sec:equilibria}, we describe the equilibria as well as their stability. 
Sections~\ref{sec:global} and \ref{sec:regions} provide a complete description of the global dynamics, including the existence of limit cycles and polycycles. In the final section
we present two families of examples when the nonlinear part of \eqref{eqR2contracting} is of degree 2 and of degree 3.

\section{Preliminary results and notation}\label{sec:prelim}

We describe the global dynamics for \eqref{eqR2contracting} depending on whether the finite degree $n$ is even or odd. 
We use the compactification  of $\RR^2$ in Chapter 5 of Dumortier {\em et al.} \cite{DLA} to describe  the dynamics at infinity  of \eqref{eqR2contracting} and to show that it determines the dynamics in $\RR^2$. The beginning of this section is devoted to recalling the Poincar\'e compactification and establishing some convenient notation that we use throughout.

Let $\mathcal{S}^2\subset \RR^3$ be the unit sphere and identify $\RR^2$ with the plane $\{ (x,y,1) \in \RR^3: x,y \in \RR\}$. Using coordinates $(z_1,z_2,z_3)$ for $\RR^3$, define charts $U_k = \{ z \in \mathcal{S}^2 : z_k > 0\}$ and $V_k = \{ z \in \mathcal{S}^2 : z_k < 0\}$ for $k=1,2,3$. 
The local maps  corresponding to these charts are
$\phi_k(z) = -\psi_k(z) = (z_m/z_k,z_n/z_k)$ for $m<n$ and $m,n \ne k$. 
Use $(u,v)$ to denote the value of the image under any of the local maps, so that the meaning of $(u,v)$ has to be determined in connection to each local chart. 
A point with coordinates $(u,v)$, $u\ne 0$ in $U_i$ or $V_i$ corresponds to the point with coordinates $(\tilde{u},\tilde{v})=\left({1}/{u},{v}/{u}\right)$ in $U_j$ or $V_j$ with $i\ne j$.
The plane $\RR^2$ is identified with the open northern hemisphere, the  Poincar\'e disk is defined as its closure.
By making $v=0$ in $U_1$, $U_2$, $V_1$, $V_2$
 we obtain the equator  $\mathcal{S}^1$ of the sphere, the {\em  circle at infinity} of the Poincar\'e disk.

A direct application of the calculations in Dumortier {\em et al.} \cite{DLA} shows that the
dynamics of  \eqref{eqR2contracting} in the Poincar\'e compactification is given,  in the  chart $U_1$, by:
\begin{equation} \label{eqPXU1}
\left\{\begin{array}{lcl}
\dot{u} &=& F(u)\\
{} &{} & \\
\dot{v} &=& -\lambda v^{n} -vQ_1(1,u)\ ,
\end{array}  \right.
\qquad\mbox{where}\qquad
F(u)=Q_2(1,u)-uQ_1(1,u) 
\end{equation}
and in the  chart $U_2$, by:

\begin{equation} \label{eqPXU2}
\left\{\begin{array}{lcl}
\dot{u} &=& G(u)\\
{} & {} & \\
\dot{v} &=&  -\lambda v^{n} -vQ_2(u,1)\ ,
\end{array}  \right. 
\qquad\mbox{where}\qquad
G(u)= Q_1(u,1)-uQ_2(u,1) .
\end{equation}
The expressions of the Poincar\'e compactification in the charts $V_1$ and $V_2$ are obtained from those in the charts $U_1$ and  $U_2$, respectively, by multiplication by $(-1)^{n-1}$.

The dynamics  at infinity of \eqref{eqR2contracting} is thus given by the restriction of either \eqref{eqPXU1} or \eqref{eqPXU2} to the flow-invariant line $(u,0)$, since the second equation is trivially satisfied for $v=0$. An equilibrium at infinity of \eqref{eqR2contracting}  is an equilibrium  $(u,0)\in \mathcal{S}^1$ of either \eqref{eqPXU1} or \eqref{eqPXU2}.
We refer to it as an {\em infinite equilibrium}, by opposition
to {\em finite equilibria} $(u,v)$, $v\ne 0$.
We refer to periodic trajectories and limit cycles as finite or infinite in the same spirit.

It is clear that the  dynamics of the restriction of either \eqref{eqPXU1} or \eqref{eqPXU2} to the flow-invariant circle at infinity $(u,0)$  does not depend on $\lambda$ and therefore, does not depend on the linear part of \eqref{eqR2contracting}. Hence, it is equivalent to 
$$
\left\{\begin{array}{lcl}
\dot{x} &=& Q_1(x,y)  \\
\mbox{} & \mbox{} & \\
\dot{y} &=& Q_2(x,y) .
\end{array}  \right.
$$
In the special case where the polynomials $Q_1$ and $Q_2$ in \eqref{eqR2contracting} have no common factor
it was established in \cite[Theorem 4.4]{AnnaL1990} that the dynamics at infinity is determined by $Q=(Q_1,Q_2)^T$. In our case the result holds without the common factor assumption.

A useful alternative description of the dynamics
 can be obtained from considering  the representation of \eqref{eqR2contracting} in polar coordinates $(x,y)=(r\cos\theta,r\sin\theta)$, $(r,\theta)\in \mathbb{R}^+\times \mathcal{S}^1$. This is given by
\begin{equation}\label{eq:polar_n}
\left\{ \begin{array}{l}
\dot{r} = \lambda r + f(\theta) r^n \\
\dot{\theta} = g(\theta) r^{n-1} \ ,
\end{array} \right. 
\end{equation}
where
\begin{equation}\label{f_polar}
f(\theta)= \left.(x,y)\cdot Q\xy \right|_{(\cos\theta,\sin\theta)} =
\cos\theta\ Q_1(\cos\theta,\sin\theta) +\sin\theta\ Q_2(\cos\theta,\sin\theta)
\end{equation}
and
\begin{equation}\label{g_polar}
g(\theta) = \left.(-y,x)\cdot Q\xy \right|_{(\cos\theta,\sin\theta)}=
 \cos\theta\ Q_2(\cos\theta,\sin\theta) -\sin\theta\ Q_1(\cos\theta,\sin\theta).
\end{equation}
Observe that $f(\theta+\pi)=(-1)^{n+1}f(\theta)$ and $g(\theta+\pi)=(-1)^{n+1}g(\theta)$.

We refer to the half-line $\theta=\theta_0$, $r>0$ as the {\em radius  $\theta=\theta_0$} and to the line $\theta=\theta_0$, $r\in\RR$ as the {\em diameter $\theta=\theta_0$}.

\begin{lemma} \label{lem:factor f}
The polynomials $Q_1$, $Q_2$ have a common linear factor if and only if $f(\theta)$ and $g(\theta)$ have a common zero.
\end{lemma}

\begin{proof}
Assume $g(\theta_0) = 0$ for some $\theta_0$, with $\cos{\theta_0} \neq 0$ (otherwise, proceed analogously for 
$\sin\theta_0\neq 0$).
Then, from \eqref{g_polar}, $g(\theta_0) = 0 \Leftrightarrow Q_2(\cos{\theta_0},\sin{\theta_0}) = \tan{\theta_0} Q_1(\cos{\theta_0},\sin{\theta_0})$. Replacing in \eqref{f_polar}, we obtain $f(\theta_0)=0 \Leftrightarrow Q_1(\cos{\theta_0},\sin{\theta_0}) = 0$. Hence, for the same $\theta_0$ it is $Q_1(\cos{\theta_0},\sin{\theta_0}) = Q_2(\cos{\theta_0},\sin{\theta_0})$ and this is true if and only if the polynomials have  a common linear factor.
\end{proof}
Note that the components of \eqref{eqR2contracting} never have common factors even when $Q_1$ and $Q_2$ do.

In order to analyse the behaviour around the circle at infinity we 
 change coordinates in \eqref{eq:polar_n} by $R=1/r$ and multiply the result  by $R^{n-1}$  obtaining the equivalent equations 
\begin{equation}\label{eq:polarinfinito}
\left\{\begin{array}{lcl}
\dot R&=&-\lambda R^n-R f(\theta)\\
\dot\theta&=&g(\theta) .
\end{array}\right.
\end{equation}
Thus, the dynamics of the vector field~\eqref{eq:polar_n} at the equator of the Poincar\'e disk is described by $\dot\theta=g(\theta)$.

\section{Equilibria}\label{sec:equilibria}
This section is concerned with the number and stability of equilibria of  \eqref{eqR2contracting}  in the Poincar\'e compactification, both finite and infinite.
We start with the behaviour at infinity.

One special case where $Q_1$ and $Q_2$ have a common factor is described in the next result.
It also implies that generically there will be either no equilibria at infinity or a finite number of them.

\begin{lemma} \label{lem:infEqinf}
There are infinitely many equilibria of \eqref{eqR2contracting}   at infinity if and only if $yQ_1(x,y)=xQ_2(x,y)$.
Moreover, in this case either there are also infinitely many finite equilibria or the only finite equilibrium is the origin.
\end{lemma}

\begin{proof}
Equilibria at infinity are points $(u,0)$ that satisfy either $F(u)=0$ or $G(u)=0$. 
Since both $F(u)$ and $G(u)$ are polynomials, in order to have infinitely many roots we must have either $F(u)\equiv 0$ or $G(u)\equiv 0$.
Direct substitution in the expressions for $F(u)$ and $G(u)$ shows that if $yQ_1(x,y)=xQ_2(x,y)$ then $F(u)=G(u)\equiv 0$. To show the converse, write

$$
Q_1(x,y)=\sum_{k=0}^n  c_{k}x^{n-k}y^k
\quad\mbox{and}\quad
Q_2(x,y)=\sum_{k=0}^n  d_{k}x^{n-k}y^k
$$
to obtain
\begin{equation}\label{eqF}
F(u)=\sum_{k=0}^n  d_{k}u^k-\sum_{k=1}^{n+1}  c_{k-1}u^k
\qquad \text{and}\qquad
G(u)=\sum_{\ell=0}^n  c_{n-\ell}u^\ell-\sum_{\ell=1}^{n+1}  d_{n-\ell+1}u^\ell .
\end{equation}
Hence, $F(u)\equiv 0$ if and only if
$$
c_{n}=0=d_{0}
\quad\mbox{and}\quad
d_{k}=c_{ k-1}\quad 
k=1,\ldots,n.
$$
Exactly the same conditions hold for $G(u)\equiv 0$.

Therefore, we can write $Q_1(x,y)=xp(x,y)$ and $Q_2(x,y)=yp(x,y)$, where  
$$
p(x,y)=\sum_{k=0}^{n-1}  c_{k}x^{n-k-1}y^k
$$ 
is a homogeneous polynomial. 
That is, $yQ_1(x,y)=xQ_2(x,y)$.
In this case, finite equilibria of \eqref{eqR2contracting} satisfy:
$$
\dot x=0\quad\Leftrightarrow\quad x=0\mbox{ or } p(x,y)=-\lambda
\qquad \mbox{and}\qquad
\dot y=0\quad\Leftrightarrow\quad y=0\mbox{ or } p(x,y)=-\lambda.
$$ 
Since the  equation $p(x,y)=-\lambda\ne 0$ has either infinitely many solutions or none, the result follows.
\end{proof}

\begin{lemma}\label{item:nparInfinito}
If $n$ is even then there is at least one pair of infinite equilibria of \eqref{eqR2contracting}.
\end{lemma}

\begin{proof}
Using \eqref{eqF}, we see that  either  the degree of $F$ is $n+1$, in which case there is at least one infinite equilibrium, or the degree is less than $n+1$.
In the second case  then $c_n=0$ and  then $G(0)=0$ and again there is an infinite equilibrium. 
  When $n$ is even, since $g(\theta)=g(\theta+\pi)$, the infinite equilibria occur in pairs.
\end{proof}

The next results provide the possible configuration of  finite equilibria.
They are similar in nature to \cite[Lemma 3.1]{CL1989}
but describe cases not covered there.

\begin{proposition}\label{lem:equilibria on radii}
There is an infinite equilibrium of \eqref{eqR2contracting} on the radius  $\theta=\theta_0$ if and only if $g(\theta_0)=0$, and in this case
\begin{enumerate}
\renewcommand{\theenumi}{(\alph{enumi})}
\renewcommand{\labelenumi}{{\theenumi}}
\item\label{invariantRadius}
the diameter $\theta=\theta_0$ is flow-invariant;
\item\label{lambdafNeg}
if $\lambda f(\theta_0)<0$ there is a unique finite equilibrium on the radius  $\theta=\theta_0$;
\item \label{lambdafPos}
if $\lambda f(\theta_0)\ge 0$  there are no finite equilibria on the radius  $\theta=\theta_0$.
\end{enumerate}
\end{proposition}

\begin{proof}
The existence of the infinite equilibrium follows from equation \eqref{eq:polarinfinito} and the invariance of the diameter is immediate from \eqref{eq:polar_n}.
Finite equilibria $(r_0,\theta_0)$ with $r_0>0$ satisfy  
$$
\left.\dot{r}\right|_{\theta_0}=\lambda r_0+f(\theta_0)r_0^n=0.
$$
Hence, $0<r_0^{n-1}=-\lambda/f(\theta_0)$, establishing \ref{lambdafNeg} and \ref{lambdafPos}.
\end{proof}

From $f(\theta +\pi)=(-1)^{n+1}f(\theta)$ and $g(\theta +\pi)=(-1)^{n+1}g(\theta)$, it follows that 
\begin{itemize}
\item  if $n$ is even, the inequalities in \ref{lambdafNeg} and \ref{lambdafPos} are reversed for $\theta_0+\pi$;
\item  if $n $ is odd \ref{lambdafNeg} and \ref{lambdafPos} also hold for $\theta_0+\pi$.
\end{itemize}
The following corollary is then immediate.

\begin{corollary}\label{cor:equilibria on radii}
Let 
$\theta_0$ be such that $g(\theta_0)=0$ in \eqref{g_polar}. Then
\begin{enumerate}
\renewcommand{\theenumi}{(\alph{enumi})}
\renewcommand{\labelenumi}{{\theenumi}}
\item 
there are two infinite equilibria on the diameter
$\theta=\theta_0$;
 \item 
 if $n$ is even then
\begin{enumerate}
\renewcommand{\theenumii}{(\roman{enumii})}
\renewcommand{\labelenumii}{{\theenumii}}
\item 
if $f(\theta_0)\ne 0$  then one of the  radii $\theta=\theta_0$ and $\theta=\theta_0+\pi$ contains a single finite equilibrium and there are no finite equilibria on other radius;
\item \label{item:radiiBii}
if $f(\theta_0)=0$ there are no finite equilibria on any of the radii $\theta=\theta_0$ and $\theta=\theta_0+\pi$;
\end{enumerate}

 \item
if $n $ is odd then
\begin{enumerate}
\renewcommand{\theenumii}{(\roman{enumii})}
\renewcommand{\labelenumii}{{\theenumii}}
\item 
if $\lambda f(\theta_0)<0$ there exists a unique finite equilibrium on each one of the radii $\theta=\theta_0$ and $\theta=\theta_0+\pi$;
\item 
if $\lambda f(\theta_0)\ge 0$ there are no finite equilibria on any of the radii $\theta=\theta_0$ and $\theta=\theta_0+\pi$.
\end{enumerate}
\end{enumerate}
\end{corollary}

 Thus, for each finite equilibrium there is a corresponding equilibrium at infinity. 
 The converse may not be true: the set of infinite equilibria may even be a continuum, with only one finite equilibrium,  as in Lemma~\ref{lem:infEqinf}.
 
 By Bezout's Theorem, if \eqref{eqR2contracting} has finitely many equilibria, then they are at most $n^2$. The next result shows that this estimate may be improved for $n>3$ (the estimate is the same if $n\leq 3$).

\begin{theorem}\label{ThCountsGeneralEquilibria}
If \eqref{eqR2contracting} has finitely many  infinite equilibria, then
\begin{enumerate}
\renewcommand{\theenumi}{(\alph{enumi})}
\renewcommand{\labelenumi}{{\theenumi}}
\item\label{teo:countsInfinito}
for all $n>1$ there are at most  $2(n+1)$ infinite equilibria;
  \item \label{teo:counts_equilibra grau n impar} 
   if $n$ is odd then  the number of  finite  equilibria away from the origin is at most $2(n+1)$;
  \item \label{teo:counts_equilibra grau n par} 
if $n$ is even then  the number of  finite equilibria away from the origin is at most $n+1$.
\end{enumerate}
\end{theorem}

\begin{proof}
We start by estimating the number of infinite equilibria.
Recall that $Q$ in \eqref{eqR2contracting} is a homogeneous polynomial vector field of degree $n$.
As in the proof of Lemma~\ref{lem:infEqinf} there are infinitely many equilibria at infinity if and only if either $F(u)\equiv 0$ or $G(u)\equiv 0$, we suppose this is not the case.
As $F(u)$ is a polynomial of degree at most $n+1$ the maximum number of infinite equilibria in the chart $U_1$ is $n+1$. Since the same holds in the chart $V_1$, the maximum number of infinite equilibria arising from the zeros of $F$ is $2(n+1)$.
The zeros of $G(u)$ yield the same infinite equilibria except for $u=0$, which does not belong to $U_1 \cup V_1$. If $G(0)=0$, by \eqref{eqF}, we have $c_n=0$ and an additional infinite equilibrium in each of the charts $U_2$ and $V_2$. However, when $c_n=0$ the degree of $F$ is $n$, rather than $n+1$, producing the same maximum number of $2(n+1)$, establishing \ref{teo:countsInfinito}.

Equilibria at infinity correspond to roots of $g(\theta)=0$. 
If the infinite equilibrium on the radius $\theta=\theta_0$ lies in the neighbourhood $U_i$ then the equilibrium on the radius $\theta=\theta_0+\pi$ lies in $V_i$.

If $n$ is odd then Proposition~\ref{lem:equilibria on radii} shows that for each one of the infinite equilibria on the radii $\theta=\theta_0$ and $\theta=\theta_0+\pi$,  there may be at most one finite equilibrium point on each one of these radii.
Therefore the total number of finite equilibria  away from the origin is at most $2(n+1)$, as in \ref{teo:counts_equilibra grau n impar}.

For even $n$,  by Proposition~\ref{lem:equilibria on radii} only one of the radii $\theta=\theta_0$ and $\theta=\theta_0+\pi$ may contain a finite equilibrium, hence the maximum number is $n+1$, establishing  \ref{teo:counts_equilibra grau n par}. 
 \end{proof}

\begin{corollary}\label{res:CountsGeneralEquilibria}
Let $Q$ be a homogeneous polynomial  vector field of degree $n=2 m+1$ and suppose \eqref{eqR2contracting} has finitely many equilibria at infinity.
If in the restriction to the circle at infinity all equilibria are either attracting or repelling then  the number of equilibria at infinity  is  a multiple of 4.
 \end{corollary}

In particular Corollary~\ref{res:CountsGeneralEquilibria} implies that when all infinite equilibria  are hyperbolic their total number is a multiple of 4.

\begin{proof}
As in the proof of Theorem~\ref{ThCountsGeneralEquilibria} we count the equilibria in the chart $U_1$ ($V_1$, respectively)  in the restriction to the circle at infinity.
Equilibria at infinity that are either attracting or repelling correspond to roots where 
$F(u)$ (respectively $G(u)$) changes sign, hence they are roots of odd multiplicity. 
If $F(u)$ has degree $2 m+2$ then there must be an even number of these roots.
In this case $G(0)\ne 0$ and  all the infinite equilibria lie in the charts $U_1$ and $V_1$.
Therefore the total number of infinite equilibria is a multiple of 4.

If the degree of $F(u)$ is $2m+2-p$ with $p\ge 1$, then from the expressions \eqref{eqF}  we get that
$$
F(u)=d_0+\sum_{k=1}^n  (d_{k}-c_{k-1})u^k-c_nu^n\ ;
$$
hence
$c_n=0$ with $d_{2m+2-p}\ne c_{2m+1-p}$ and $ c_{k-1}=d_k$ for $k\ge 2m+3-p$.
Then
$$
G(u)=c_n+\sum_{\ell=1}^n  (c_{n-\ell}-d_{n-\ell+1})u^\ell-  d_{0}u^{n+1}
=(c_{2m+1-p}-d_{2m+2-p})u^p+\sum_{\ell=p+1}^n  (c_{n-\ell}-d_{n-\ell+1})u^\ell-  d_{0}u^{n+1},
$$
since if $k=n-\ell+1$ then $c_{n-\ell}=c_{k-1}=d_k=d_{n-\ell+1}$.
Therefore $G(0)=0$ with multiplicity $p$.

If $p$ is even, then $G$ does not change sign at $u=0$, contradicting the hypothesis.
If $p$ is odd, then the number of  roots of $F(u)$ is odd, say $2k+1$, so there are $4k+2$ infinite equilibria in the charts $U_1$ and $V_1$.
In this case $G$  changes sign at $u=0$, corresponding to two equilibria not in the charts $U_1$ and $V_1$,
and the total number of infinite equilibria is again a multiple of 4.
\end{proof}

The next result concerns the stability of the equilibria of  \eqref{eqR2contracting}. 
In the case of infinite equilibria it repeats results in  Proposition 4.1 (c) in \cite{AnnaL1990}, which we include here for ease of reference.

\begin{proposition}\label{prop:stability}
If  \eqref{eqR2contracting} has finitely many equilibria at infinity, then:
\begin{enumerate}
\renewcommand{\theenumi}{(\alph{enumi})}
\renewcommand{\labelenumi}{{\theenumi}}
 \item \label{stab:realEigs}
the linearisation of \eqref{eqPXU1} and \eqref{eqPXU2} at any (finite or infinite) equilibrium has real eigenvalues,
in particular, all equilibria are either topological nodes (attractors or repellors), saddles or saddle-nodes;
\item \label{stab:radialInfinite}
the infinite equilibrium on the radius $\theta=\theta_0$ is radially attracting if  either $ f(\theta_0)> 0$ or $ f(\theta_0)= 0$ and $\lambda>0$;  it is
radially repelling if  either $ f(\theta_0)< 0$ or if $ f(\theta_0)= 0$ and $\lambda<0$;
 \item \label{stab:all}
the stability in  the angular direction of an equilibrium (finite or infinite) on the radius $\theta=\theta_0$ is determined by the sign of $g'(\theta_0)$, and in the radial direction, for finite equilibria, by the sign of $-\lambda$;
\item  \label{stab:saddleNode}
an equilibrium $(u_0,v_0)$ (finite or infinite) in the chart $U_1$ (respectively $U_2$) is a saddle-node if and only if $u_0$ is a root of even multiplicity of $F(u)$ in \eqref{eqPXU1}  (respectively $G(u)$ in \eqref{eqPXU2}).
\end{enumerate}
\end{proposition}

\begin{proof}

In both \eqref{eqPXU1} and \eqref{eqPXU2} the expression for $\dot u$ does not depend on $v$, hence the Jacobian matrix is triangular and the eigenvalues are real, proving \ref{stab:realEigs}.

Using polar coordinates \eqref{eq:polar_n} at a finite equilibrium $(r_0,\theta_0)$ we have 
$0<r_0^{n-1}=-\lambda / f(\theta_0)$
and the Jacobian matrix of \eqref{eq:polar_n} is the triangular matrix
\begin{equation}\label{JacobianoPolar}
J(r_0,\theta_0)=
\left(
  \begin{array}{cc}
    -\lambda (n-1) & f'(\theta_0)r_0^{n} \\
    0 & g'(\theta_0) \left(\frac{-\lambda}{f(\theta_0)}\right) \\
  \end{array}
\right).
\end{equation}
Hence,  the radial direction is an eigenspace corresponding to the non-zero eigenvalue $-\lambda (n-1)$, whose sign is the opposite of that of $\lambda$, proving  \ref{stab:radialInfinite} and the second part of \ref{stab:all}.

The stability of a finite equilibrium $(r_0,\theta_0)$ in the angular direction depends on the sign of $g'(\theta_0)$ since  $0<r_0^{n-1}=-\lambda / f(\theta_0)$. It follows from equation \eqref{eq:polarinfinito}  that this stability coincides with that at infinity, proving \ref{stab:all}.

If  $f(\theta_0)\ne 0$ then for small $R$ the sign of $\dot R$ in equation \eqref{eq:polarinfinito} is the same as the sign of $-f(\theta_0)$.
When $f(\theta_0)= 0$ then  $\dot R$ has the  sign of $-\lambda$  for small $R>0$,
proving \ref{stab:radialInfinite}. Note that for $f(\theta)=0$ the infinite equilibrium is not hyperbolic.

For an infinite equilibrium, statement \ref{stab:saddleNode} is proved  in \cite[Proposition 4.1 (c)]{AnnaL1990}.
At a finite equilibrium the radial eigenvalue of $J(r_0,\theta_0)$ is never zero, and \ref{stab:saddleNode} holds for the angular direction, by \ref{stab:all}.
\end{proof}

\section{Stability of the origin}\label{sec:global}
 The aim of this section and the next one  is to describe how the dynamics of \eqref{eqR2contracting} at infinity constrains the geometry of flow-invariant sets other than equilibria. 
 We start by the cases when the dynamics is simpler.
The next result is a synthesis of Theorem 2 in  \cite{BLS2013} and Theorem A in \cite{CGP1997}, as well as Theorems 2 and 3  in \cite{GYZ}, we include it here for ease of reference.

\begin{theorem}[\cite{BLS2013,CGP1997,GYZ}]\label{teo:Bende}
For \eqref{eqR2contracting}  with $n>1$ and $\lambda\ne 0$:
\begin{enumerate}
\renewcommand{\theenumi}{(\roman{enumi})}
\renewcommand{\labelenumi}{{\theenumi}}
\item\label{item:BendeAtMostOne}
there is at most one non-constant  finite periodic solution and its trajectory surrounds the origin;
\item\label{item:BendeA}
if $n$ is even then there are
no finite periodic trajectories surrounding the origin;
\item\label{item:BendeB}
if $n$ is odd and $g(\theta)=0$ for some $\theta\in[0,2\pi)$ then there are
 no  finite periodic trajectories surrounding the origin;
 \item\label{item:BendeC}
if $n$ is odd and $g(\theta)\ne 0$ for all $\theta\in[0,2\pi)$ then the origin is the only equilibrium of \eqref{eqR2contracting} and there is a  finite limit cycle surrounding the origin if and only if
$\lambda {\mathcal I}=\lambda\int_0^{2\pi}\frac{f(\theta)}{|g(\theta)|}d\theta<0$;
\item\label{item:BendeImplicit}
if $n$ is odd and $g(\theta)\ne 0$ for all $\theta\in[0,2\pi)$, then the limit cycle at infinity is stable if and only if ${\mathcal I}>0$;
\item\label{item:hyperbolic}
if there is a limit cycle it is hyperbolic.
\end{enumerate}
\end{theorem}

From this result, Proposition~\ref{lem:equilibria on radii} and Corollary~\ref{cor:equilibria on radii}  we obtain  necessary conditions for trivial dynamics.

\begin{corollary}\label{coro:globalattractor}
If the origin is a globally attracting ($\lambda < 0$) or  repelling ($\lambda > 0$) equilibrium point  of \eqref{eqR2contracting} then one of the following conditions holds:
\begin{enumerate}
\renewcommand{\theenumi}{(\alph{enumi})}
\renewcommand{\labelenumi}{{\theenumi}}
\item\label{item:necess f>0 impar}
$n$ is odd with  $g(\theta)\ne 0$ for all $\theta\in[0,2\pi)$ and $\lambda\int_0^{2\pi}\frac{f(\theta)}{|g(\theta)|}d\theta\ge 0$;
\item\label{item:necess impar}
$n$ is odd and  $\lambda f(\theta)\geq 0$  
whenever $g(\theta)=0$;
\item\label{item:necess par}
 $n$ is even and  $f(\theta)=0$ whenever $g(\theta)=0$.
\end{enumerate}
\end{corollary}

\begin{proof}
By Proposition~\ref{lem:equilibria on radii}, Corollary~\ref{cor:equilibria on radii} and Theorem~\ref{teo:Bende}  the conditions on $f(\theta)$ and $g(\theta)$ of \ref{item:necess f>0 impar}--\ref{item:necess par} are equivalent to the origin being the only finite equilibrium of \eqref{eqR2contracting} and the non-existence of finite limit cycles.
Hence the conditions are necessary. 
\end{proof}

\begin{proposition}\label{prop:globalattractor}
The origin is a globally attracting ($\lambda < 0$) or  repelling ($\lambda > 0$) equilibrium point  of \eqref{eqR2contracting} if one of the following conditions holds:
\begin{enumerate}
\renewcommand{\theenumi}{(\alph{enumi})}
\renewcommand{\labelenumi}{{\theenumi}}
\item\label{item:attractor f>0 impar}
$n$ is odd with  $g(\theta)\ne 0$ for all $\theta\in[0,2\pi)$ and $\lambda\int_0^{2\pi}\frac{f(\theta)}{|g(\theta)|}d\theta\ge 0$;
\item\label{item:attractor impar}
$n$ is odd, $Q_1$ and $Q_2$ have no common linear factor, and $\lambda f(\theta)\geq 0$ whenever $g(\theta)=0$.
\end{enumerate}
\end{proposition}

\begin{proof}
We have already established  that these conditions imply that there are no finite equilibria except for the origin.
Condition~ \ref{item:attractor f>0 impar} implies that there are no infinite equilibria.
 To show the condition is sufficient we establish the absence of non-trivial finite periodic trajectories. 
By Theorem~\ref{teo:Bende}~\ref{item:BendeA} if such trajectories existed they would surround the origin.
Under condition \ref{item:attractor f>0 impar}, Theorem~\ref{teo:Bende}~\ref{item:BendeC} and \ref{item:BendeImplicit} guarantee no such trajectory exists and the stability of the limit cycle at infinity is the opposite of that of the origin.

For condition \ref{item:attractor impar}, Theorem~\ref{teo:Bende}~\ref{item:BendeB} and \ref{item:BendeC} imply that there are no finite limit cycles. 
Since $Q_1$ and $Q_2$ have no common factor,  Lemma~\ref{lem:factor f} implies that $ f(\theta_0)\ne 0$ when $g(\theta_0)=0$, hence $\lambda f(\theta_0)>0$ at these points.
In order to analyse the behaviour around the equilibria at infinity we use equation \eqref{eq:polarinfinito} around $R=0$.
On the restriction to the invariant ray $\theta=\theta_0$  the infinite equilibrium $(R,\theta)=(0,\theta_0)$ is hyperbolic and repelling for $\lambda<0$, so no finite trajectory can have it as accumulation point and the result follows by the Poincar\'e-Bendixson Theorem.
For $\lambda>0$ the same holds with reversed time.
\end{proof}

If $n$ is even then by Corollary~\ref{cor:equilibria on radii}  the origin is the only finite equilibrium of \eqref{eqR2contracting} if and only if $f(\theta)=0$ whenever $g(\theta)=0$.
By Theorem~\ref{teo:Bende} \ref{item:BendeA} there are also no finite periodic trajectories.
 The problem is making sure there are no  finite trajectories connecting the infinite equilibria.

The next example, constructed by Fabio Scalco Dias, illustrates the need for the  condition that $Q_1$ and $Q_2$ have no common linear factor in Proposition~\ref{prop:globalattractor}  \ref{item:attractor impar}.

\begin{figure}[hht]
 \begin{center}
  \includegraphics[width=0.2\linewidth]{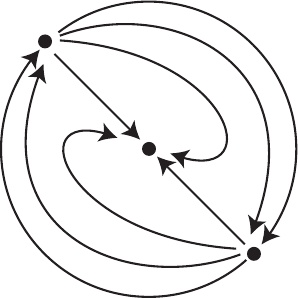}
  \end{center}
  \caption{ Dynamics of Example~\ref{ex:grau3} on the Poincar\'e disk.}
  \label{fig:examples}
\end{figure}

\begin{Ex}[F. Scalco Dias]\label{ex:grau3}
Let
\begin{equation}\label{eq:exgrau3}
\left\{\begin{array}{lcl}
\dot x&=& -x +(x^3+y^3)\\
\dot y&=&-y -(x^3+y^3)
\end{array}\right.
\qquad n=3\quad \lambda=-1.
\end{equation}
Then $Q_2\xy=-Q_1\xy$ with $Q_1\xy=x^3+y^3=\dfrac{1}{4}(x+y)\left[3(x-y)^2+(x+y)^2\right]$, so 
$x+y$ is the only common linear factor.
Since 
$$
g(\theta) 
=\left.-\dfrac{1}{4}(x+y)^2\left(3(x-y)^2+(x+y)^2\right)\right|_{(\cos\theta,\sin\theta)}\le 0
$$
then $g(\theta) =0$ if and only if $\theta=3\pi/4$ or $\theta=-\pi/4$, otherwise $g(\theta)<0$.
At these points $f(\theta)=0$
since
$$
f(\theta)=\left.\dfrac{1}{4}(x-y)(x+y)\left(3(x-y)^2+(x+y)^2\right) \right|_{(\cos\theta,\sin\theta)}
$$
 so on the corresponding radius $\dot r=-r<0$.

To see the behaviour at infinity we use equation \eqref{eq:polarinfinito}.
On the line $R=0$   the infinite equilibrium $(R,\theta)=(0,-\pi/4)$ is a saddle-node attracting on the $\theta>-\pi/4$ side.
On the ray $R>0$, $\theta=-\pi/4$ this equilibrium is a non hyperbolic source $\dot R =R^3$.
However, in the region $R^2<f(\theta)$, close to the line at $\infty$, $\dot R<0$ and this is where a trajectory
connecting the two infinite equilibria
$(0,3\pi/4)$ and $(0,-\pi/4)]$ may exist.

Let  $F(R,\va)$ be the vector field associated to \eqref{eq:polarinfinito} in the new coordinate $\va=\theta+\pi/4$.
To see if there is such a connection we blow up the vector field $F(R,\va)$ around $(R,\va)=(0,0)$ in the direction $\va$.
To simplify the calculations we use the degree 4 Taylor expansions of $f$ and $g$, that are:
\begin{equation}\label{eq:taylor}
f(\va)=3\va-4\va^3+O(\va^5)
\qquad \text{and} \qquad
g(\va)=-3\va^2+3\va^4+O(\va^6) .
\end{equation}
Then 
the rescaled directional blow up (see \cite{DLA}) is 
$\tilde{F}(R,\va)=\dfrac{1}{R^2}(DM)^{-1}\cdot F\circ M(R,\va)$, 
where $M(R,\va)=(R,R\va)$. 
Using the Taylor expansions \eqref{eq:taylor}, it is given by:
$$
\tilde{F}(R,\va)=\begin{pmatrix}
R-3\va+4R^2\va^3\\
-\va-R\va^4
\end{pmatrix}+{\mathcal O}(R^4\va^5)\ ,
$$
hence
$$
D\tilde{F}(R,\va)=\begin{pmatrix}
1+8\va^3&-3+12R^2\va^2\\
-\va^4&-1-4R\va^3
\end{pmatrix}+{\mathcal O}(R^3\va^4)
\quad\Rightarrow\quad
D\tilde{F}(0,0)=\begin{pmatrix}
1&-3\\ 0&-1
\end{pmatrix}
$$
with eigenvalues: $1$ and $-1$, eigenvectors $(1,0)^T$ and $(3,2)^T$ respectively.

Therefore there are trajectories of \eqref{eq:polarinfinito} in the first quadrant $R>0$, $\theta>-\pi/4$ going into the infinite equilibrium $(0,-\pi/4)$ and the origin is not a global attractor, as shown in Figure~\ref{fig:examples}.

Note that this example is very degenerate, since the linear change of coordinates $(u,v)=\left(x+y,x-y\right)$
brings it to the form
$$
\left\{\begin{array}{lcl}
\dot u&=&-u+0\\
\dot v&=&-v+2(x^3+y^3) \ .
\end{array}\right.
$$
\end{Ex}

 Examples of degree 2 with  similar dynamics may be found in the book \cite{Artes}.

  \section{Invariant regions and polycycles}\label{sec:regions}
In  this section we discuss the dynamics of \eqref{eqR2contracting} when $g(\theta_0)=0$ for some $\theta_0\in[0,2\pi)$.
This condition, by Proposition~\ref{lem:equilibria on radii}, implies the existence of two infinite equilibria.
It also implies the flow-invariance of the diameter $\theta=\theta_0$. 
 Suppose  $g(\theta)$ is not identically zero  and let $\ba_1$ and $\ba_2$ be two consecutive infinite equilibria corresponding to two consecutive zeros of $g(\theta)$, $\theta_1$ and $\theta_2$. These define a flow-invariant set given in polar coordinates by $\left\{(r,\theta)\ : \ r\in\RR,\ \theta_1\le\theta\le \theta_2\right\} $, 
with $0<\theta_2-\theta_1\le \pi$.
We refer to this flow-invariant set as a {\em flow-invariant cone}  (or simply a {\em cone}) and to its non-negative ($r\ge 0$) component as a  {\em half-cone}.

Therefore, to equilibria at infinity there corresponds a division of the plane in invariant half-cones between  consecutive invariant radii.
An upper bound for the number of half-cones can be obtained from Theorem~\ref{ThCountsGeneralEquilibria}.

\begin{corollary}\label{cor:half-cones}
If the homogeneous polynomial $Q$ has degree $n>1$ then there are at most $2(n+1)$ half cones that are invariant under the flow of  \eqref{eqR2contracting}.
If $n$ is even there is at least one pair of invariant half-cones.
\end{corollary}

\begin{proof}
Each pair of infinite equilibria, satisfying $g(\theta_0)=0=g(\theta_0+\pi)$ determines a flow-invariant diameter.
From the proof of Theorem~\ref{ThCountsGeneralEquilibria} it follows that there are at most $2(n+1)$ such pairs, independently of the parity of the degree $n$.
If $n$ is even then by Proposition~\ref{lem:equilibria on radii} there is at least one pair of infinite equilibria, giving rise to one pair of invariant half-cones.
\end{proof}

\begin{figure}[hht]
  \centering
  \includegraphics[width=.8\linewidth]{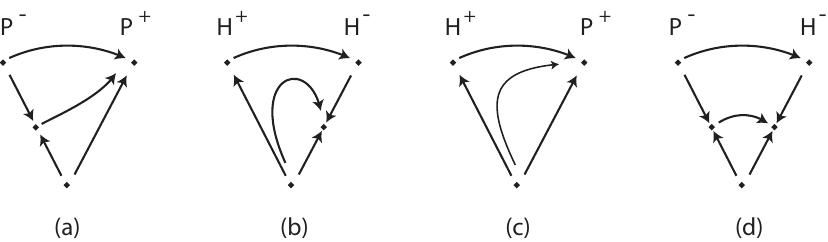}
  \caption{Dynamics in invariant half-cones in Theorem~\ref{teo:NewGlobalDynamics} for  $\lambda>0$. 
 In case \ref{teoGlobalPP} the unstable manifold of the finite equilibrium splits the half-cone in two basins of attraction.  In case \ref{teoGlobalPH} there is a robust heteroclinic connection between two finite equilibria. 
In the remaining cases the origin is a global repellor for the interior of the half-cone.
  }\label{fig: invariantregions}
\end{figure}

An  infinite equilibrium  $\ba$  and its associated invariant radius determine two angles in a neighbourhood in the half-plane $v\ge 0$.
Denote by $P^-$ an angle that is repelling in both the line at infinity and in the radius, $P^{+}$ an angle that is  attracting in both the line at infinity and in the radius, $H^{-}$ an  angle that is  attracting in  the line at infinity and radially repelling and $H^{+}$ an angle that is repelling in  the line at infinity and radially attracting, 
see Figure~\ref{fig: invariantregions}.

A {\em parabolic sector} is a subset of the   Poincar\'e disk bounded by either a $P^+$ or a $P^-$ angle such that the dynamics is qualitatively equivalent to that of a node.
A {\em hyperbolic sector} is bounded  by either a $H^+$ or a $H^-$ angle with dynamics qualitatively equivalent to a saddle.
A region bounded by an angle may contain more than one sector as in Figure~\ref{fig:examples}.

We will  use below this classification of sectors to determine the dynamics of  \eqref{eqR2contracting} on the invariant half-cones. 

If $\ba$ is a hyperbolic equilibrium of \eqref{eq:polarinfinito}, i.e. if $Q_1$ and $Q_2$ have no common linear factor and $g'(\ba)\ne 0$,  then the $P^\pm$ are parabolic sectors and the $H^\pm$ are hyperbolic sectors.

\begin{lemma}\label{lema:cones}
If $g(\theta)$ is not identically zero the possible pairs of angles at two consecutive infinite equilibria are:
\begin{enumerate}
\renewcommand{\theenumi}{(\alph{enumi})}
\renewcommand{\labelenumi}{{\theenumi}}
\item \label{teoGlobalPP}
$P^-$ and $P^+$;
\item \label{teoGlobalHH}
$H^+$ and $H^-$;
\item \label{teoGlobalHP}
$H^+$  and $P^+$;
\item\label{teoGlobalPH}
$P^-$ and $H^-$.
\end{enumerate}
\end{lemma}

\begin{proof}
Since the infinite equilibria are consecutive, the flow at infinity goes from one of the infinite equilibria to the other. Then, if the  angles are of the same type, $H$ or $P$, they cannot have the same sign.
If the angles are of different types then, of course, one is $P^\pm$.
If it is $P^-$ then the other   angle, which is  $H^\pm$,  is attracting in the angular component and must therefore be radially repelling. That is, it is $H^-$. 
 Analogously, if  one angle  is $P^+$ then the other   angle is repelling in the angular component. It must therefore be radially attracting, that is, it is $H^+$.
\end{proof}

\begin{figure}[hht]
  \centering
   \includegraphics[width=.8\linewidth]{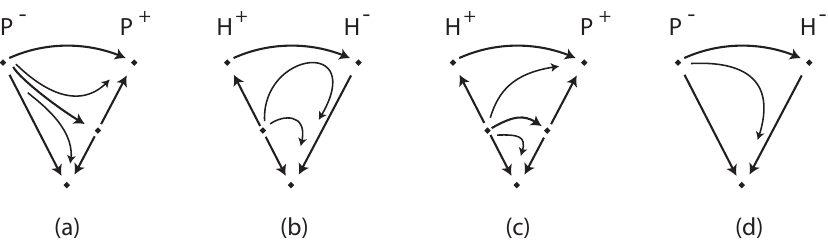}
  \caption{Dynamics in invariant half-cones in Theorem~\ref{teo:NewGlobalDynamics} for  $\lambda<0$. In case \ref{teoGlobalPP} the stable manifold of the finite equilibrium splits the half-cone in two basins of attraction.  In case \ref{teoGlobalHP} there is a robust heteroclinic connection between two finite equilibria.  In the remaining cases the origin is a global attractor for the interior of the half-cone.}\label{fig:halfconesLambdaNegative}
\end{figure}

\begin{theorem} \label{teo:NewGlobalDynamics}
If all the infinite equilibria of the vector field \eqref{eqR2contracting} are hyperbolic, then the dynamics near them determines the global dynamics on the plane.
For a trajectory at infinity connecting two consecutive infinite equilibria $\ba_1$ and $\ba_2$ the dynamics in the half-cone $C$ determined by them is described in Figures~\ref{fig: invariantregions} and \ref{fig:halfconesLambdaNegative}, where  the possible pairs of  angles at $\ba_1$ and $\ba_2$ are those in Lemma~\ref{lema:cones}.
\end{theorem}

\begin{proof}
Assume that $\lambda>0$ so that the origin is repelling as in Figure~\ref{fig: invariantregions}.
If  all the infinite equilibria are hyperbolic, then equation \eqref{eq:polarinfinito} implies that $g(\theta)$ is not identically zero, and there are finitely many infinite equilibria.
Because of Proposition~\ref{lem:equilibria on radii}, one is the maximum number of finite equilibria on the radius $\theta = \theta_i$ corresponding to $\ba_i$, $i=1,2$. We have that
\begin{itemize}
\item   an angle delimited by a repelling radius exists at $\ba_i$ if and only if exactly one finite equilibrium exists on the radius $\theta = \theta_i$;
\item    an angle delimited by an attracting radius exists at $\ba_i$ if and only if there are no finite equilibria on the radius $\theta = \theta_i$.
\end{itemize}
Since 
to a finite equilibrium there corresponds an infinite one, and because
$\ba_1$ and $\ba_2$ are consecutive infinite equilibria, it follows
 that there are no finite equilibria in the interior of $C$.

The admissible dynamics are depicted in Figure~\ref{fig: invariantregions}. 
\smallskip

For $\lambda <0$ the origin is attracting and the dynamics in $C$ may be obtained by changing time as $s=-t$. The phase portraits are depicted in Figure~\ref{fig:halfconesLambdaNegative}.
\end{proof}

It remains to see what happens if the infinite equilibrium $\ba$ at $\theta_0$ is not hyperbolic.
This may happen if either $f(\theta_0)=0$ (i.e. $Q_1$ and $Q_2$ have a common linear factor) or if $g'(\theta_0)=0$.
In this case the angle formed by the line at infinity and the radius may comprise more than one dynamically defined sector.
We start by establishing the dynamics around an angle of type $P$.

\begin{lemma}\label{lema:parabolicsector}
If $\ba$ is an infinite equilibrium of the vector field \eqref{eqR2contracting}  that determines an angle that is  attracting (respectively repelling) in both the line at infinity and in the radius, then the angle is a parabolic sector.
\end{lemma}

\begin{proof}
We show that  the non hyperbolic equilibrium $\ba$ with coordinates $(R,\theta)=(0,\theta_0)$ attracts (respectively repels) all points in the intersection of the half-cone with  an open rectangle around $\ba$. 
Assume that $\lambda>0$ (for $\lambda<0$ reverse time), and for simplicity, that the sector is locally $\theta>\theta_0$.

First note that if $f(\theta)$ is identically zero, then the first equation in \eqref{eq:polarinfinito} reduces to $\dot R=-\lambda R^n<0$  for $R>0$.
In particular this holds for $\theta_0$ which implies that $\ba$ is radially attracting.
It follows  by hypothesis that the angle is $P^+$ and hence $\ba$ is attracting in the line at infinity, therefore $g(\theta)<0$ for $\theta>\theta_0$ close to $\theta_0$.
From \eqref{eq:polarinfinito} and since there is no other equilibrium around $\ba$,   it follows that  $\ba$ attracts all trajectories with $R>0$ and $\theta>0$ close to $\theta_0$, so the angle is a parabolic  sector.

If $f(\theta_0)\ge 0$ then at the radius at $\theta_0$ we have $\dot R=-\lambda R^n-Rf(\theta_0)<0$ for small $R>0$, so $\ba$ is radially attracting and the angle is of type $P^+$ with $g(\theta)<0$ for $\theta>\theta_0$ close to $\theta_0$.
By continuity, there are $R_1>0$ and $\theta_1>\theta_0$ such that if $0<R\le R_1$ and $\theta_0<\theta\le \theta_1$ we have $\dot R<0$ and $\dot\theta=g(\theta)<0$.
Therefore the rectangle $(R,\theta)\in [0,R_1]\times[\theta_0,\theta_1]$ is positively flow-invariant so it must contain the $\omega$-limit set of all trajectories in it.
Since the only equilibrium in the rectangle is $\ba$,  this must be the $\omega$-limit set,  so the angle is a parabolic    sector.
The same reasoning may be applied if $f(\theta_0)< 0$ with $\dot R=-\lambda R^n-Rf(\theta_0)<0$ for small $R>0$.

If $f(\theta_0)< 0$ with $\dot R=-\lambda R^n-Rf(\theta_0)>0$ for small $R>0$ then $\ba$ is radially repelling and the angle is of type $P^-$ with $g(\theta)>0$ for $\theta>\theta_0$ close to $\theta_0$.
Then there are $R_1>0$ and $\theta_1>\theta_0$ such that if $0<R\le R_1$ and $\theta_0<\theta\le \theta_1$ we have $\dot R>0$ and $\dot\theta=g(\theta)>0$.
 Therefore the rectangle $(R,\theta)\in [0,R_1]\times[\theta_0,\theta_1]$ is negatively flow-invariant so it must contain the $\alpha$-limit set of all trajectories in it.
Since the only equilibrium in the rectangle is $\ba$,  this must be the $\alpha$-limit set,  so the angle is a parabolic sector.
\end{proof}

An angle of type $H$ may contain more than one sector, as we saw in Example~\ref{ex:grau3}.
Figure~\ref{fig:examples} shows two infinite equilibria such that the angle of type $H^-$ contains one hyperbolic and one parabolic sector.
The next result establishes that the only possibility is the situation of that example.

\begin{lemma}\label{lema:hyperbolicsector}
If $\ba$ is an infinite equilibrium of the vector field \eqref{eqR2contracting}  that determines an angle that is  attracting  (respectively repelling) in  the line at infinity and radially repelling  (respectively attracting), then the angle contains a hyperbolic sector and at most one parabolic sector.
The last case only occurs if $f(\theta_0)=0$ where $(0,\theta_0)$ are the coordinates $(R,\theta)$ at $\ba$.
\end{lemma}

\begin{proof}
Let $(R,\theta)=(0,\theta_0)$ be the coordinates of the non hyperbolic equilibrium $\ba$  and without loss of generality consider an angle defined locally by $\theta>\theta_0$.
Again we deal with the case $\lambda>0$, for $\lambda <0$ the result follows by reversing time.

If $f(\theta)$ is identically zero, then, as in the proof of Lemma~\ref{lema:parabolicsector} it follows that $\ba$ is radially attracting, hence the angle is of type $H^+$, and there exists $\theta_1>\theta_0$ such that $g(\theta)>0$ for $\theta_0<\theta\le\theta_1$. For some $R_1>0$, consider the  rectangle 
$(R,\theta)\in[0,R_1]\times[\theta_0,\theta_1]$. 
Its boundary consists of the flow-invariant segments $(0,\theta)$ and $(R,\theta_0)$, of the segment $(R_1,\theta)$ where the vector field points in  (since $\dot R=-\lambda R^n<0$) and of the segment $(R,\theta_1)$ where the vector field points out. 
Trajectories starting at $\{R_1\}\times(\theta_0,\theta_1)$ go into the rectangle.
The only equilibrium in the rectangle is  $\ba=(0,\theta_0)$ and that point cannot be their $\omega$-limit, since $\dot \theta=g(\theta)>0$, so the trajectories must go out through the segment $[0,R_1]\times\{\theta_1\}$.
Therefore, the angle is a hyperbolic sector.

If $f(\theta_0)>0$, since 
at the radius at $\theta_0$ we have $\dot R=-\lambda R^n-Rf(\theta_0)<0$,  then  $\ba$ is radially attracting and the angle is of type $H^+$. 
By continuity, there is   $\theta_1>\theta_0$ such that $f(\theta)>0$ and $g(\theta)<0$ for $\theta\in[\theta_0,\theta_1]$, hence, for any $R_1>0$, in the rectangle $(R,\theta)\in[0,R_1]\times[\theta_0,\theta_1]$ we have
$\dot R<0$ and $\dot \theta>0$. 
The same arguments used in the case $f(\theta)\equiv 0$ show that the  angle is a hyperbolic sector.

If $f(\theta_0)<0$ then,   for $R<R_1=\left(-f(\theta_0)/\lambda \right)^{1/(n-1)}$, at the radius at $\theta_0$ we have  $\dot R=-\lambda R^n-Rf(\theta_0)>0$ hence $\ba$  is radially repelling and 
by hypothesis it is attracting in the line at infinity, therefore the angle is of type $H^-$.
%
%
By continuity, there is a rectangle $(R,\theta)\in[0,R_2]\times[\theta_0,\theta_1]$  with $R_2<R_1$, where $\dot R>0$ and $\dot\theta=g(\theta)<0$. 
Its boundary consists of the flow-invariant segments $(0,\theta)$ and $(R,\theta_0)$, of the segment $(R_2,\theta)$ where the vector field points out  (since $\dot R>0$) and of the segment $(R,\theta_1)$ where the vector field points in. 
Trajectories starting at $(R,\theta_1)$ go into the rectangle.
The only equilibrium in the rectangle is  $\ba=(0,\theta_0)$ and that point cannot be their $\omega$-limit, since $\dot R>0$, so the trajectories must go out through the segment $(R_2,\theta)$.
 Similarly, all trajectories starting at points in the rectangle with $R>0$ satisfy $\dot \theta<0$ so $\ba=(0,\theta_0)$ cannot be their $\alpha$-limit.
Therefore, the angle is a hyperbolic sector.

If $f(\theta_0)=0$ then $\ba$ is radially attracting  and the angle is of type $H^+$.
If $f(\theta)>0$ for $\theta>\theta_0$ close to $\theta_0$, then, as in the case $f(\theta_0)>0$, there is a rectangle $(R,\theta)\in[0,R_1]\times[\theta_0,\theta_1]$ where $f(\theta)>0$ and $g(\theta)>0$,  hence $\dot R<0$ and $\dot \theta>0$ and the  angle is a hyperbolic sector, as in the case $f(\theta_0)>0$.

If $f(\theta_0)=0$ with $f(\theta)<0$ for $\theta>\theta_0$ close to $\theta_0$, then $\dot R=0$ on the curve $R=R(\theta)=\left(-f(\theta)/\lambda \right)^{1/(n-1)}$, with $\dot R>0$ for $R<R(\theta)$.
Since for some $\theta_1>\theta_0$ we have $g(\theta)>0$ for $\theta_0<\theta\le\theta_1$,
then, for trajectories starting at $(R,\theta_1)$ with $0<R<R(\theta_1)$, the coordinates $R(t)$ and $\theta(t)$ are monotonically increasing functions.
If the trajectory  meets
the curve $R=R(\theta)$ it crosses it from the region where $\dot R<0$ from where it must arrive at the rectangle $(R,\theta)\in[0,R(\theta_1)]\times[\theta_0,\theta_1]$ through the side $(R(\theta_1),\theta)$.
However, there may be trajectories that never cross the curve $R(\theta)$ and therefore must have $\ba$ as $\alpha$-limit.
In this case the angle at $\ba$ consists of two sectors, one hyperbolic $H^+$ and the other parabolic $P^-$.
Since $\dot\theta>0$ at all points in the interior of the rectangle,  then no trajectory can have $\ba$ as $\omega$-limit.
This precludes the existence of an {\em elliptic sector}, filled with loops having both  $\alpha$- and $\omega$-limit $\ba$.
\end{proof}

Note that  in Lemmas~\ref{lema:parabolicsector} and \ref{lema:hyperbolicsector} the fact that  $g'(\theta_0)=0$  is irrelevant to the dynamics on the angle. 

\begin{figure}[hht]
  \centering
   \includegraphics[width=.8\linewidth]{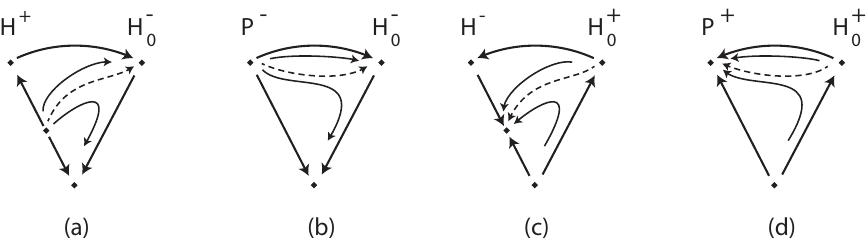}
  \caption{Dynamics in invariant half-cones in Theorem~\ref{teo:globalDynamics}  when there is an angle $H_0^\pm$ consisting of two sectors. The dashed line is the trajectory separating the sectors, it  splits the half-cone in two basins of attraction or repulsion. In (a) and (b)  $\lambda<0$, in (c) and (d)  $\lambda>0$.  In cases (b) and (d) there is a  parabolic sector  at infinity connecting  the two infinite equilibria as in Example~\ref{ex:grau3}.  }\label{fig:halfconesNonHyperbolic}
\end{figure}

\begin{theorem} \label{teo:globalDynamics}
If  the vector field \eqref{eqR2contracting}  has finitely many equilibria  at infinity, then 
the dynamics near them determines its global dynamics on the plane.

For a trajectory at infinity connecting two consecutive infinite equilibria $\ba_1$ and $\ba_2$ the dynamics in the half-cone $C$ determined by them is described in Figures~\ref{fig: invariantregions}, \ref{fig:halfconesLambdaNegative} and \ref{fig:halfconesNonHyperbolic}.
\end{theorem}

\begin{proof}
Assume that $\lambda>0$, so the origin is repelling and consider a half-cone $C$ defined by two consecutive infinite equilibria $\ba_1$ and $\ba_2$ at the angles $\theta_1<\theta_2$, respectively.
If both $\ba_1$ and $\ba_2$ are hyperbolic, then by Theorem~\ref{teo:NewGlobalDynamics} the dynamics around them determines the dynamics in $C$.
Suppose that at least one of the equilibria is not hyperbolic.
We claim that the admissible dynamics in the half-cone is  one of those shown either in Figure~\ref{fig:halfconesLambdaNegative}  or in Figure~\ref{fig:halfconesNonHyperbolic}.

To establish the claim, note that in the  proof of Theorem~\ref{teo:NewGlobalDynamics} the hypothesis that the infinite equilibria are hyperbolic is only used to ensure that all the angles in the half-cone consist of a single sector, either parabolic or hyperbolic.
By Lemmas~\ref{lema:parabolicsector} and \ref{lema:hyperbolicsector} this is true of almost all half-cones,
so  the arguments in that proof can be applied except in the case where one of the angles contains two sectors.
In this case from the proof of Lemma~\ref{lema:hyperbolicsector} it follows that the angle is of type $H^+$, we denote it by $H_0^+$.
Without loss of generality, let  $\ba_1$ be the infinite equilibrium at this angle, so  $f(\theta_1)=0$.
By Lemma~\ref{lema:cones} the  angle defined by $\ba_2$ is either $H^-$ or $P^+$ and hence it consists of a single sector. From now on we suppose that $H_0^+$ contains a repelling parabolic sector.

If the angle at $\ba_2$ is $H^-$, then the radius at $\theta_2$ contains a single finite equilibrium.
This equilibrium is the only possibility in $C$ for the $\omega$-limit set of the  trajectories in the repelling sector in $H_0^+$ since there are no finite equilibria in the interior of $C$ nor  in the radius at $\theta_1$ and the origin is repelling. Thus, the dynamics is that of  Figure~\ref{fig:halfconesNonHyperbolic} (c).

If the angle at  $\ba_1$ is $P^+$ then none of the radii at $\theta_1$ or $\theta_2$ contains a finite equilibrium.
The only possibility in $C$ for the $\omega$-limit set of the  trajectories in the repelling sector in $H_0^+$ is the attractor $\ba_1$ and the dynamics  is that of  Figure~\ref{fig:halfconesNonHyperbolic} (d).
\end{proof}

 A {\em polycycle} is a flow-invariant simple connected curve in the plane containing at least one equilibrium point and not going through the origin. In the special case when all the trajectories that are not equilibria have the same orientation it is called a {\em heteroclinic cycle}.

When one sector is parabolic and the other hyperbolic a finite equilibrium may exist in each radius. If these exist, their radial stability is opposite to that of the origin. Thus the connection between these two finite equilibria is robust as it is either of saddle-sink ($\lambda>0$) or saddle-source ($\lambda <0$) type. See Figures~\ref{fig: invariantregions} \ref{teoGlobalPH} and \ref{fig:halfconesLambdaNegative} \ref{teoGlobalHP}. 
Note that for a polycycle to exist, all sectors must  be  alternatingly  either $H^-$ and $P^-$
if $\lambda >0$ or  $H^+$ and $P^+$ if $\lambda <0$.
The next results provide more detail concerning polyclycles and heteroclinic cycles.

\begin{proposition} \label{prop:poly-cycle}
The dynamics of \eqref{eqR2contracting} exhibits a polycycle if and only if $n$ is odd, there is at least one pair of infinite equilibria and $\lambda f(\theta)<0$ whenever $g(\theta)=0$. Moreover, the polycycle is globally attracting 
(respectively repelling) if $\lambda >0$ (respectively $\lambda <0$).
\end{proposition}

\begin{proof}
By Lemmas~\ref{lema:parabolicsector} and \ref{lema:hyperbolicsector} if  $\lambda f(\theta)<0$ when $g(\theta)=0$ then any angle at infinity defined by an infinite equilibrium consists of a single sector, either parabolic or hyperbolic.

  If $n$ is odd and $\lambda f(\theta)<0$ when $g(\theta)=0$  it follows by Proposition~\ref{lem:equilibria on radii} that each invariant radius contains a finite equilibrium. 
  On each invariant half-cone one of the sectors is $P^s$ and the other is $H^s$ where $s$ is the sign of $-\lambda$.
  By Theorem~\ref{teo:globalDynamics} there is a trajectory in the half-cone connecting the two finite equilibria, forming a polycycle. Its dynamics in the radial direction is given by the sign of $-\lambda$. 
  Thus, the polycycle is globally attracting (respectively repelling) if $\lambda >0$ (respectively $\lambda <0$).
  
Conversely, if $n$ is even, then for each $\theta_0$ such that $g(\theta_0)=0$, by Proposition~\ref{lem:equilibria on radii}, there is always a radius where there are no finite equilibria.
Hence there is a pair of consecutive finite equilibria for which there is no connecting trajectory, as it would have to cross an invariant radius, and there is no polycycle.

 The same argument shows that if $n$ is odd and there is a polycycle, then there must be a finite equilibrium on each invariant radius.
Proposition~\ref{lem:equilibria on radii} implies that $\lambda f(\theta)<0$ whenever $g(\theta)=0$.
\end{proof}

Another way of stating Proposition~\ref{prop:poly-cycle} is that any (finite) polycycle is a copy of the polycycle at infinity, with the radial stability inverted.
All the connections of a polycycle for \eqref{eqR2contracting}  are robust, by \ref{teoGlobalHP} and \ref{teoGlobalPH} of Lemma~\ref{lema:cones}.

\begin{corollary}\label{cor:heteocycle} 
 A polycycle of \eqref{eqR2contracting} is a heteroclinic cycle if and only if all infinite equilibria are local minima or local maxima of $g$.
\end{corollary}

\begin{proof}
In a heteroclinic cycle all the equilibria must be connected by trajectories with the same orientation.
All the equilibria must be saddle-nodes.
Therefore the sign of $g(\theta)$ must be always the same. 
\end{proof}

A simple interpretation of Corollary~\ref{cor:heteocycle} follows from Figures~\ref{fig: invariantregions} and \ref{fig:halfconesLambdaNegative}.
In order  to have a heteroclinic cycle all half-cones must be as in Figure~\ref{fig: invariantregions} (d) for $\lambda>0$ (respectively Figure~\ref{fig:halfconesLambdaNegative} (c) for $\lambda<0$) and hence  all connections must be $P^-\seta H^-$ if $\lambda>0$ (respectively $H^+\seta P^+$ if $\lambda<0$).
This is the only result in this section where the derivative $g'(\theta_0)$ is relevant.

In parametrised systems a limit cycle may be created through a saddle-node bifurcation of the equilibria in a heteroclinic cycle.
The next result shows that indeed this happens within the class we are studying.

\begin{proposition}\label{prop:bifurca}
A 1-parameter perturbation of a heteroclinic cycle for \eqref{eqR2contracting} creates a limit cycle.
\end{proposition}

\begin{proof}
Let $Q_\varepsilon(x,y) = \left(Q_1(x,y)+\varepsilon P_1(x,y), Q_2(x,y)+\varepsilon P_2(x,y)\right)$ 
be homogeneous of degree $n$ 
with $( P_1(x,y),P_2(x,y)) = (-y^n, x^n)$.
Define $f_\varepsilon(\theta)$ and $g_\varepsilon(\theta)$ as in \eqref{f_polar} and \eqref{g_polar}, respectively. We have
\begin{eqnarray*}
f_\varepsilon(\theta) & = & f(\theta) + \varepsilon \sin\theta \cos\theta (\cos^{n-1}\theta -\sin^{n-1}\theta )\\
g_\varepsilon(\theta) & = & g(\theta) + \varepsilon (\cos^{n+1}\theta +\sin^{n+1}\theta ).
\end{eqnarray*}
Assume that a heteroclinic cycle exists for the dynamics of \eqref{eqR2contracting} when $\varepsilon = 0$. By Proposition~\ref{prop:poly-cycle}  and Corollary~\ref{cor:heteocycle} it must be that the degree of $Q$ is odd, say $n=2m+1$, and $g(\theta)$ has constant sign. Without loss of generality, we assume in this proof that $g(\theta)\ge 0$ for all $\theta$ and that $\lambda>0$, the other cases being analogous. Hence, we also have $f(\theta)<0$ when $g(\theta)=0$.

Note that when $n=2m+1$, for $\varepsilon>0$ we have
 $g_\varepsilon(\theta) >0$ for all $\theta$. Hence, the vector field defined by $Q_\varepsilon$ has no infinite equilibria and only the origin is a finite equilibrium.
The proof is completed by showing that
\begin{equation}\label{eq:integral}
\dpt \int_0^{2\pi} \dfrac{f_\varepsilon (\theta)}{|g_\varepsilon (\theta)|} d\theta < 0.
\end{equation}
Then the hypotheses of Theorem~\ref{teo:Bende} \ref{item:BendeC} are satisfied and for small $\varepsilon >0$ a limit cycle exists.
Choose $\varepsilon>0$ small enough so that 
$f_\varepsilon(\theta)<0$ for $\theta$ in intervals containing the $\theta_i$ where $g(\theta_i)=0$.
When $\varepsilon\to 0$ we have $g_\varepsilon(\theta_i)\to 0$.
Since $\lim_{\varepsilon\to 0}1/{g_\varepsilon(\theta_i)}=+\infty$,
the contribution of these intervals to the value of the integral in \eqref{eq:integral} increases when $\varepsilon$ decreases and  the integral is negative for small $\varepsilon>0$.
\end{proof}

When  \eqref{eqR2contracting} has only one pair of infinite equilibria the saddle-node bifurcation described in Proposition~\ref{prop:bifurca} holds for an open set of homogeneous vector fields near $Q$.
This is clearly not true if there is more than one pair of equilibria, since the collapse of the saddle-nodes may not be simultaneous for a generic 1-parameter perturbation.

\section{Examples}\label{sec:examples}

The results of  the previous sections are applied here to obtain phase portraits of \eqref{eqR2contracting} for nonlinearities of degrees 2 and 3. 
The dynamics of \eqref{eqR2contracting} at infinity depends only on the nonlinear part and by
Theorem~\ref{teo:globalDynamics} it determines the global dynamics, except when there are no infinite equilibria.

A linear change of coordinates and a time rescaling that  preserves orientation transforms the equation 
$\dot X=\lambda X+Q(X)$, with $\lambda\ne 0$, where $X=(x,y)$ and $Q$ is a non-zero homogeneous  vector field, into $\dot X=\tilde\lambda X+\beta {\mathcal Q}(X)$ with $\beta=\pm1$, where $\tilde\lambda$ has the same sign as $\lambda$ and ${\mathcal Q}$ is a homogeneous polynomial of the same degree as $Q$.
This is true because the linear part $\lambda X$ of the equation commutes with every linear map of $\RR^2$.

Thus, we can apply
the results of Date \cite{Date} and of Cima and Llibre \cite{AnnaL1990} on the classification of homogeneous polynomial vector fields of 
 degrees 2 and 3, under  linear changes of coordinates and a rescaling of time.
 For higher degree the interested reader may apply the same procedure to the classification of Collins \cite{Collins}.
 
 Given  a homogeneous  vector field  $Q = (Q_1, Q_2)$  of degree $n$  on the plane, let $\GG\xy=xQ_2\xy-yQ_1\xy$, a homogeneous form of degree $n+1$.
  It follows from equation \eqref{eq:polarinfinito} that $g(\theta)=\GG (\cos\theta,\sin\theta)$ determines the dynamics of \eqref{eqR2contracting} at infinity.
For the dynamics near infinity we also need the expression of $f(\theta)=\FF(\cos\theta,\sin\theta)$ where
$\FF\xy=xQ_1\xy+yQ_2\xy$.

\subsection{Example: nonlinearities of degree 2}\label{subsec:deg2}
As an illustration of our results we obtain all the possible phase portraits for  \eqref{eqR2contracting}, when the nonlinearities are quadratic. 
This repeats  results presented in Art\'es {\sl et al.}\cite{Artes}, except for the possibility of existence of  parabolic sectors discussed in Lemma~\ref{lema:hyperbolicsector}.
We use the classification of homogeneous cubic forms $\GG \xy$  from \cite{Date} and \cite{AnnaL1990}. 
In order to apply these classifications
 we compute  in the next result the general form of the vector fields of the form \eqref{eqR2contracting} with nonlinearities of degree 2 that share the same expression for $\GG $.
 Its proof is a direct computation.
\begin{lemma}\label{lem:binary-to-vf-2}
Let $Q = (Q_1, Q_2)$ be a homogeneous quadratic vector field on the plane defining \eqref{eqR2contracting}. 
 If $\GG (x,y)=xQ_2\xy-yQ_1\xy= a_0x^3+ a_1x^2y+a_2xy^2+a_3 y^3$, then 
 there exist $q_1, q_2\in\RR$ such that \eqref{eqR2contracting} has the form: 
 $$
\left\{\begin{array}{lcl}
\dot{x} &=& \lambda x + q_1x^2+(q_2-a_2)xy-a_3y^2\\
\dot{y} &=& \lambda y + a_0x^2+(a_1+q_1)xy+q_2y^2
\end{array}  \right.
$$
with $\FF\xy= q_1x^3+(q_2+a_0-a_2)x^2y+(a_1-a_3+q_1)xy^2+q_2y^3$.
\end{lemma}

The next result recovers \cite{Date} with the normal forms for $\GG $ from 
Theorem 1.4 in \cite{AnnaL1990}.
Here $\dot X=\tilde\lambda X- {\mathcal Q}(X)$ is mapped into $\dot X=\tilde\lambda X+{\mathcal Q}(X)$ by a rotation of $\pi$ around the origin, hence we do not need to add $\beta=\pm 1$ to the canonical forms.

\begin{figure}[hhb]
\parbox{30mm}{
\begin{center}
\includegraphics[width=\linewidth]{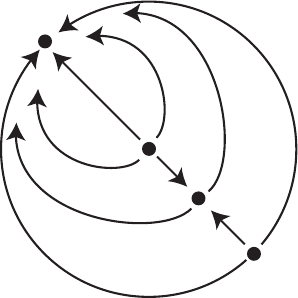}\\
 (A)
\end{center}
}
\qquad
\parbox{30mm}{
\begin{center}
\includegraphics[width=\linewidth]{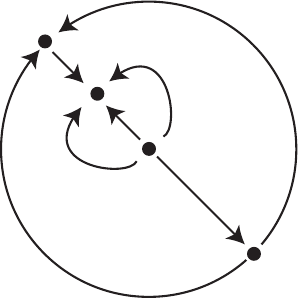}\\
 (B)
\end{center}
}
\qquad
\parbox{30mm}{
\begin{center}
\includegraphics[width=\linewidth]{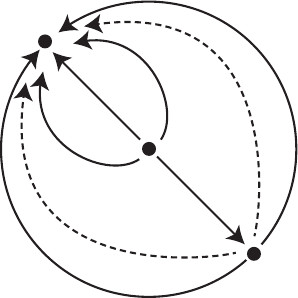}\\
 (C)
\end{center}
}
\caption{Phase portraits on the Poincar\'e disk for normal form \ref{Date1} in Proposition~\ref{prop:10vf-2} with $\lambda>0$: (A)  $q_1-q_2<1$; (B) $q_1-q_2>1$; (C) $q_1-q_2=1$.
The dashed line in the case  (C) where $q_1-q_2=1$ indicates  the possibility of
 a parabolic sector at the infinite equilibrium $\theta=-\pi/4$. 
  Whether this sector does or does not exist has to be decided by other methods, see Appendix~\ref{app:correspondence}.
 }
\label{fig:Date1}
\end{figure}

\begin{figure}[hht]
\includegraphics[width=0.65\linewidth]{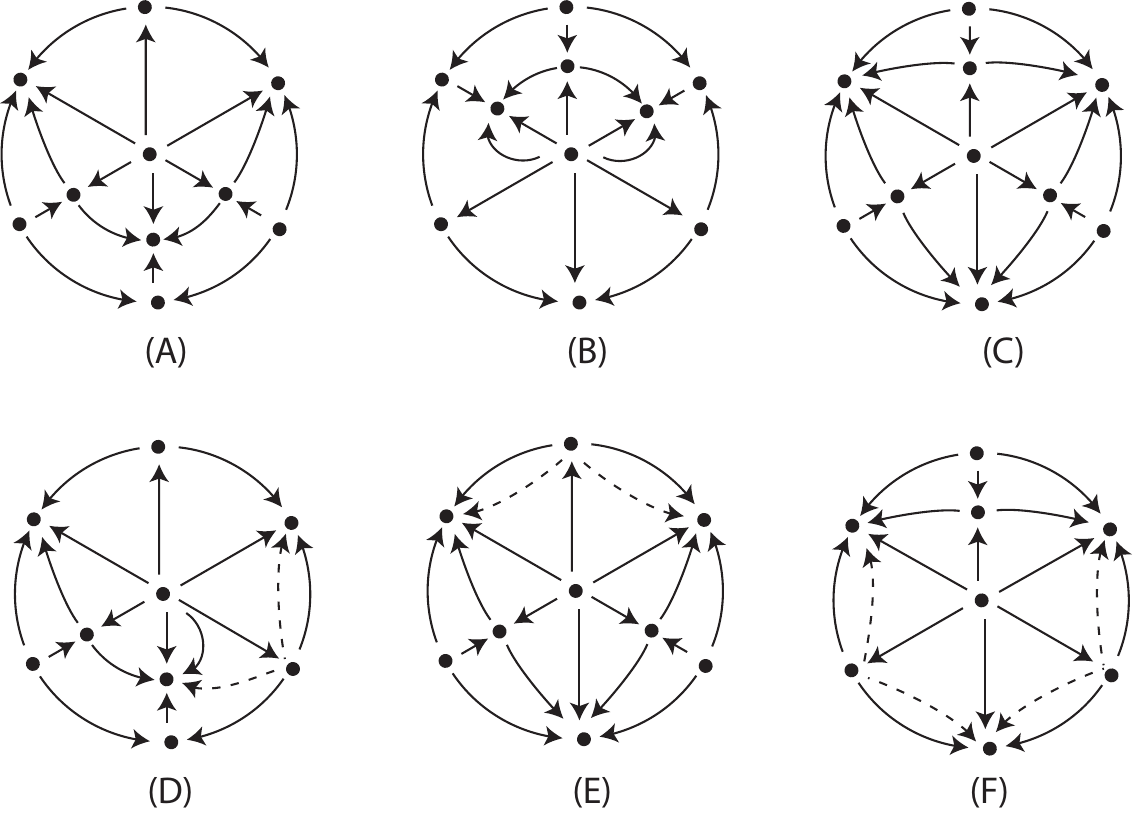}
\caption{Phase portraits on the Poincar\'e disk for normal form  \ref{Date2} in Proposition~\ref{prop:10vf-2} with 
$\lambda>0$.
Conditions on the parameters $q_1$, $q_2$ for (A)--(F) are given  in Appendix~\ref{app:deg2}.
The dashed lines in the  transition cases (D)--(F) indicate the possibility of parabolic sectors at  infinite equilibria
  that may or may not exist.
Their existence  has to be decided by other methods, see Appendix~\ref{app:correspondence}.}  
\label{fig:Date2}
\end{figure}

\begin{figure}[hht]
\includegraphics[width=0.65\linewidth]{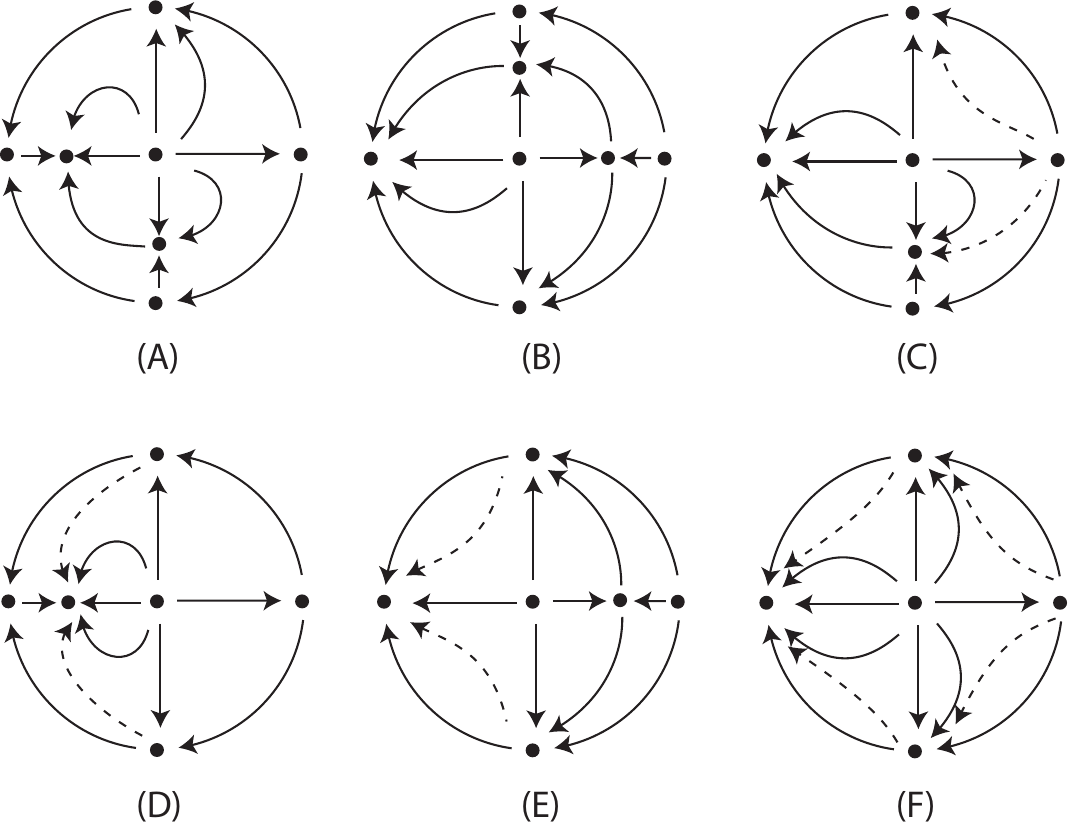}
\caption{Phase portraits on the Poincar\'e disk for normal form  \ref{Date3} in Proposition~\ref{prop:10vf-2} with $\lambda>0$.
Conditions on the parameters $q_1$, $q_2$:
(A) $q_1>0$, $q_2\ne 0$;
(B) $q_1<0$, $q_2\ne 0$;
(C) $q_1=0$, $q_2\ne 0$;
(D) $q_1>0$, $q_2=0$;
(E) $q_1<0$, $q_2= 0$;
(F) $q_1=0$, $q_2=0$.
The dashed lines in the last four transition cases indicate the possibility of parabolic sectors at  infinite equilibria  that may or may not exist.
Their existence  has to be decided by other methods, see Appendix~\ref{app:correspondence}.} 
\label{fig:Date3}
\end{figure}

\begin{figure}[hht]
\begin{center}
\includegraphics[width=0.65\linewidth]{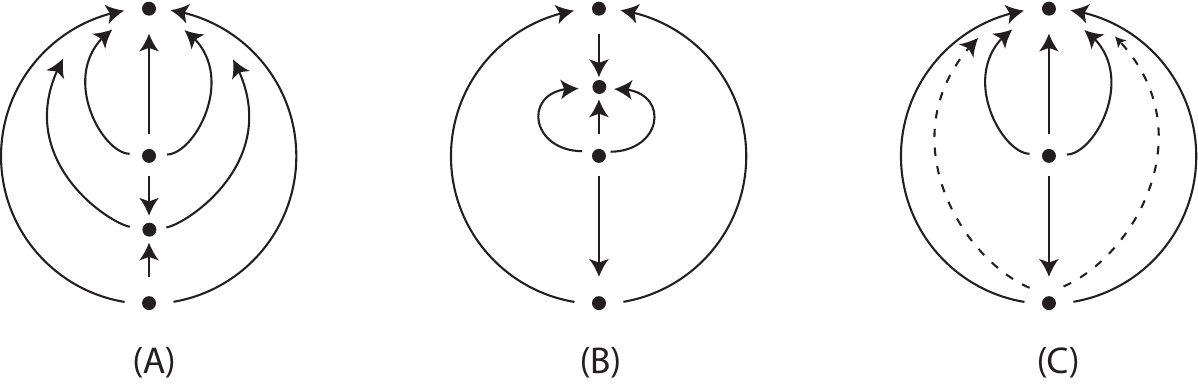}
\end{center}
\caption{Phase portraits on the Poincar\'e disk for normal form
 \ref{Date4} in Proposition~\ref{prop:10vf-2} with $\lambda>0$. 
 Conditions on the parameters $q_1$, $q_2$:
(A) $q_2>0$;
(B) $q_2<0$;
(C) $q_2=0$.
 The dashed lines in  case (C) indicate the possibility of   parabolic sectors at an infinite equilibrium  that may or may not exist.
Their existence  has to be decided by other methods, see Appendix~\ref{app:correspondence}.}  
\label{fig:Date4}
\end{figure}

\begin{figure}[hhh]
\begin{center}
\includegraphics[width=30mm]{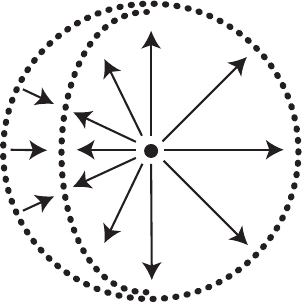}
\end{center}
\caption{Phase  portrait on the Poincar\'e disk for normal form  \ref{Date5} in Proposition~\ref{prop:10vf-2} with $\lambda>0$ and $(q_1,q_2)\ne(0,0)$. The  dotted lines are continua  of equilibria, both finite and infinite.
 The finite equilibria lie on the line $\lambda+q_1x+q_2y=0$.
  If $q_1=q_2=0$ then the equation is linear.}
\label{fig:Date345}
\end{figure}

\begin{proposition}\label{prop:10vf-2}
For each homogeneous quadratic vector field $Q = (Q_1, Q_2)$  defining \eqref{eqR2contracting} 
there exists  a linear change of coordinates and an orientation preserving reparametrisation of time that transforms \eqref{eqR2contracting} into only one of the following canonical forms, where $\tilde{\lambda}$ and $\lambda$ have the same sign:
\begin{enumerate}
\renewcommand{\theenumi}{(\roman{enumi})}
\renewcommand{\labelenumi}{{\theenumi}}
  \item\label{Date1}
  $\left\{\begin{array}{lclcl}
\dot{x} &=& \tilde{\lambda} x +q_1x^2+q_2xy-y^2 
&\qquad \qquad &\GG\xy=x^3+y^3 \\
\dot{y} &=& \tilde{\lambda} y  + x^2+q_1xy+q_2y^2 
&\qquad \qquad &\FF\xy=q_1x^3+(q_2+1)x^2y+(q_1-1)xy^2+q_2y^3;
\end{array}  \right.$
\medbreak

  \item\label{Date2}
  $\left\{\begin{array}{lclcl}
\dot{x} &=& \tilde{\lambda} x +q_1x^2+(q_2+3)xy 
&\qquad \qquad &\GG\xy=x(x^2-3y^2)\\
\dot{y} &=& \tilde{\lambda} y  + x^2+q_1xy +q_2y^2
&\qquad \qquad &\FF\xy=q_1x^3+(q_2+4)x^2y+q_1xy^2+q_2y^3;
\end{array}  \right.$
\medbreak

  \item\label{Date3}
  $\left\{\begin{array}{lclcl}
\dot{x} &=& \tilde{\lambda} x +q_1x^2+q_2xy 
&\qquad \qquad &\GG\xy=3x^2y\\
\dot{y} &=& \tilde{\lambda} y  + (q_1+3)xy+q_2y^2 
&\qquad \qquad &\FF\xy=q_1x^3+q_2x^2y+(q_1+3)xy^2+q_2y^3;
\end{array}  \right.$

\medbreak

  \item\label{Date4}
  $\left\{\begin{array}{lclcl}
\dot{x} &=& \tilde{\lambda} x+q_1x^2+q_2xy 
&\qquad \qquad &\GG\xy=x^3\\
\dot{y} &=& \tilde{\lambda} y  + x^2+q_1xy +q_2y^2
&\qquad \qquad &\FF\xy=q_1x^3+(q_2+1)x^2y+q_1xy^2+q_2y^3;
\end{array}  \right.$
\medbreak

  \item\label{Date5}
  $\left\{\begin{array}{lclcl}
\dot{x} &=& \tilde{\lambda} x+q_1x^2+q_2xy 
&\qquad \qquad &\GG\xy=0\\
\dot{y} &=& \tilde{\lambda} y +q_1xy+q_2y^2
&\qquad \qquad &\FF\xy=q_1x^3+q_2x^2y+q_1xy^2+q_2y^3 .
\end{array}  \right. $
\end{enumerate}
For $\lambda>0$ these canonical forms give rise to a minimum of 16 qualitatively different phase portraits depicted in Figures~\ref{fig:Date1}--\ref{fig:Date345}.
\end{proposition}

The phase portraits for these canonical forms may now be obtained using Proposition~\ref{lem:equilibria on radii}, its corollary, Lemma~\ref{lema:hyperbolicsector} and Theorem~\ref{teo:globalDynamics}.
They are shown in Figures~\ref{fig:Date1}--\ref{fig:Date345}
for $\lambda>0$. The phase portraits for case \ref{Date1} (Figure~\ref{fig:Date1}) are topologically equivalent to those of case \ref{Date4} (Figure~\ref{fig:Date4}), although in the latter case the infinite equilibria are not hyperbolic.
Calculations are given in Appendix~\ref{app:deg2}.
The finer classification in \cite[Theorem 7.2]{Artes} enumerates 40 possible configurations of quadratic systems with a star node. 
They determine explicitly when  there is a parabolic sector at an infinite equilibrium
 and
distinguish between parabolic sectors where nearby trajectories are tangent to the radius and those where they are tangent to the line at infinity, while our results do not account for this.
 The correspondence between their classification and ours is given in Appendix~\ref{app:correspondence}
where we use it to show that there are 17 qualitatively different phase portraits.


\newpage

\phantom{lixo}
\newpage


\subsection{Example: nonlinearities of degree 3}\label{subsec:deg3}
Cima and Llibre in \cite{AnnaL1990} classify  homogeneous binary forms $\GG \xy$ of degree 4  along with  homogeneous polynomial vector fields of  degree 3, under linear changes of coordinates and a rescaling of time.
Since the rescaling allows a reversal of time, we multiply their vector fields by a  term  $\beta=\pm 1$, except in the cases when the change of sign can be achieved through the parameters $p_1,p_2,p_3$ and $\alpha$.
We use this classification to obtain the dynamics at infinity as a starting point  for the following list of canonical forms of homogeneous cubic vector field $Q = (Q_1, Q_2)$  defining \eqref{eqR2contracting}.

\begin{figure}[hht]
  \centering
   \includegraphics[width=30mm]{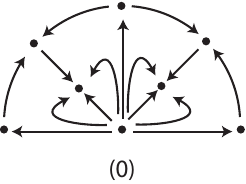}
   \qquad
  \includegraphics[width=30mm]{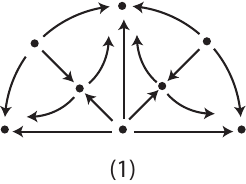}
  \qquad   
  \includegraphics[width=30mm]{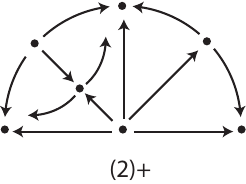}
  \\  \vskip5mm
 \includegraphics[width=30mm]{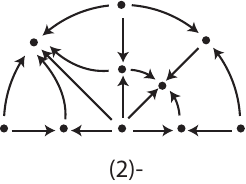}
  \qquad
 \includegraphics[width=30mm]{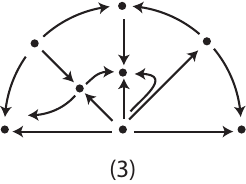}
   \qquad   
 \includegraphics[width=30mm]{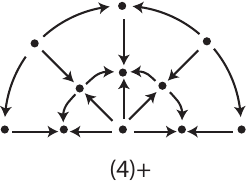}
  \\  \vskip5mm
 \includegraphics[width=30mm]{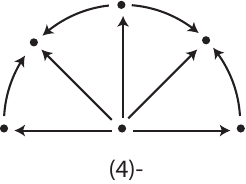}
  \qquad
  \includegraphics[width=30mm]{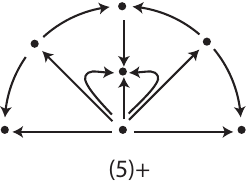}
  \qquad 
   \includegraphics[width=30mm]{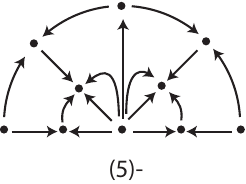}
  
  \caption{Phase portraits for normal form \ref{Ana1} in Proposition~\ref{prop:deg3} for $\tilde\lambda>0$,
  when $Q_1$ and $Q_2$ have no common factor.
 Each diagram in  \cite{AnnaL1990}[Figure 5.1]  gives rise to two cases, depending on the choice of sign, numbering refers to their cases. In cases (0),  (1) and (3) the two choices give equivalent phase portraits, case (0) is missing from  \cite{AnnaL1990}.
  Only half the disk is shown, the other half is obtained by  a rotation of $\pi$ around the origin.
  The equator of the Poincar\'e sphere is a global attractor in case (4)- and  there is a polycycle in case (4)+. 
  }\label{fig: caseI}
\end{figure}

\begin{figure}[hht]
  \centering
   \includegraphics[width=30mm]{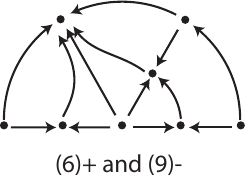}
  \qquad   
  \includegraphics[width=30mm]{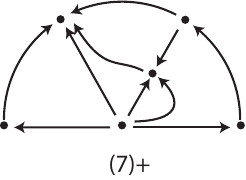}
 \qquad
 \includegraphics[width=30mm]{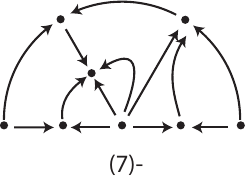}
 \\  \vskip5mm
 \includegraphics[width=30mm]{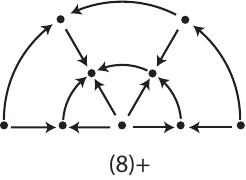}
   \qquad   
 \includegraphics[width=30mm]{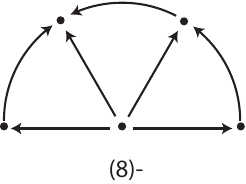}
  \qquad
 \includegraphics[width=30mm]{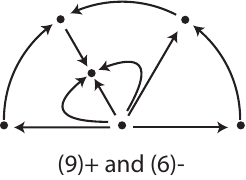}
  
  \caption{Phase portraits for normal form \ref{Ana5} in Proposition~\ref{prop:deg3} for $\tilde\lambda>0$    when $Q_1$ and $Q_2$ have no common factor.
 Each diagram in  \cite{AnnaL1990}[Figure 5.1]  gives rise to two cases, depending on the choice of sign, numbering refers to their cases. 
 Cases (6)$\pm$ and (9)$\mp$ coincide.
  The equator of the Poincar\'e sphere is a global attractor in case (8)- and  there is a polycycle in case (8)+. 
  Conventions as in Figure~\ref{fig: caseI}.
  }\label{fig: caseV}
\end{figure}

\begin{figure}[hht]
  \centering
   \includegraphics[width=30mm]{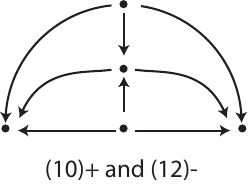}
   \qquad
  \includegraphics[width=30mm]{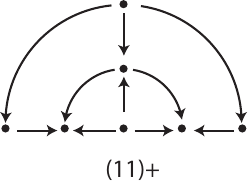}
  \qquad   
  \includegraphics[width=30mm]{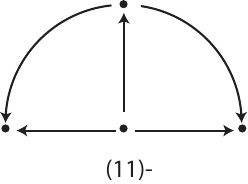}
 \qquad
 \includegraphics[width=30mm]{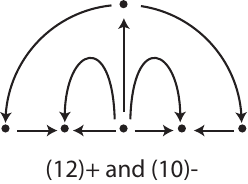}
 \\  \vskip5mm
 \includegraphics[width=30mm]{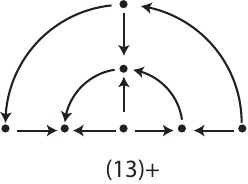}
   \qquad   
 \includegraphics[width=30mm]{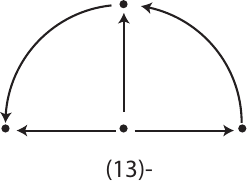}
  \qquad
 \includegraphics[width=30mm]{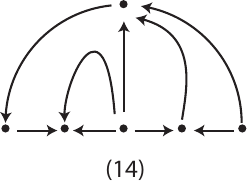}
  
  \caption{Phase portraits  with $\tilde\lambda>0$ for normal forms \ref{Ana3}, \ref{Ana7} and \ref{Ana8} in Proposition~\ref{prop:deg3}, where there are 4 infinite equilibria when $Q_1$ and $Q_2$ have no common factor.
  Numbering refers to cases  in  \cite{AnnaL1990}[Figure 5.1]. 
 Cases (10)--(12) correspond to  both \ref{Ana3} and  \ref{Ana8}  with (10)$\pm$ and (12)$\mp$ coinciding, (14)+ and (14)- are the same.
 Cases (13) and (14) arise for \ref{Ana7}.  
 The equator of the Poincar\'e sphere is a global attractor in cases (11)- and (13)-.  
 There is a polycycle in cases (11)+ and (13)+. 
 Conventions as in Figure~\ref{fig: caseI}.
  }\label{fig: caseIII}
\end{figure}

\begin{figure}[hht]
  \centering
   \includegraphics[width=30mm]{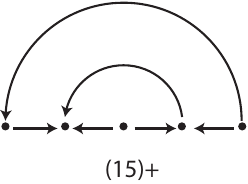}
   \qquad
  \includegraphics[width=30mm]{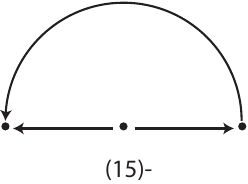}
  
  \caption{Phase portraits with $\tilde\lambda>0$ for normal forms \ref{Ana4} and \ref{Ana9}   in Proposition~\ref{prop:deg3}, with two infinite equilibria,  when $Q_1$ and $Q_2$ have no common factor.
  The equator of the Poincar\'e sphere is a global attractor in case (15)- and   there is a heteroclinic cycle  in case (15)+.
  Conventions as in Figure~\ref{fig: caseI}.
  }\label{fig:casesIVandIX}
\end{figure}

\begin{figure}[hht]
  \centering
  \includegraphics[width=30mm]{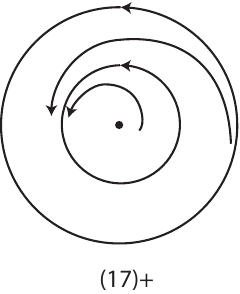}
   \qquad   
 \includegraphics[width=30mm]{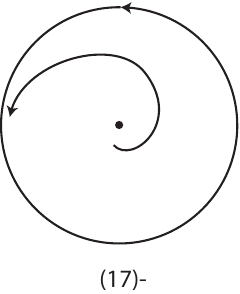}
  
  \caption{Phase portraits for normal forms \ref{Ana2} and \ref{Ana6}, with no infinite equilibria,  in Proposition~\ref{prop:deg3} for $\tilde\lambda>0$.
  The equator of the Poincar\'e sphere is a global attractor in case (17)- and  there is a  limit cycle in case (17)+.
  Conventions as in Figure~\ref{fig: caseI}.
  }\label{fig: casesIiandVI}
\end{figure}

\begin{figure}
 \centering
 \parbox{31mm}{\centering
   \includegraphics[width=30mm]{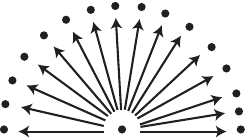}\\
   \ref{CXD+T+}  and   \ref{CXD0T+} 
}
  \qquad   
 \parbox{31mm}{\centering
   \includegraphics[width=30mm]{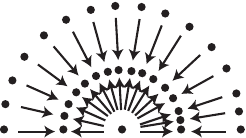}\\
   \ref{CXD+T-}
}
  \qquad 
   \parbox{31mm}{\centering
   \includegraphics[width=30mm]{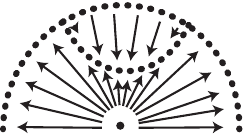}\\
   \ref{CXD-}
}
  \qquad 
   \parbox{31mm}{\centering
   \includegraphics[width=30mm]{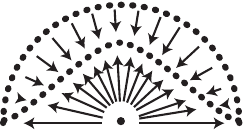}\\
   \ref{CXD0T-}
}
\caption{Phase portraits  with $\tilde\lambda>0$ for normal form \ref{Ana10} in Proposition~\ref{prop:deg3}, when all points at infinity are equilibria.  Cases  \ref{CXD+T+}--\ref{CXD0T-} as in Lemma~\ref{lem:caseX}.
 Conventions as in Figure~\ref{fig: caseI}.
  }\label{fig: caseX}
\end{figure}

\begin{proposition}\label{prop:deg3} 
For each non-zero homogeneous cubic vector field $Q = (Q_1, Q_2)$  defining \eqref{eqR2contracting} there exists  a linear change of coordinates and an orientation preserving reparametrisation of time that transforms \eqref{eqR2contracting} into only one of the following canonical forms, where $\tilde{\lambda}$ and $\lambda$ have the same sign:

\begin{enumerate}
\renewcommand{\theenumi}{(\Roman{enumi})}
\renewcommand{\labelenumi}{{\theenumi}}
  \item\label{Ana1}
  $\left\{\begin{array}{lcl}
\dot{x} &=& \tilde{\lambda} x +\beta \left[x (p_1x^2+p_3y^2)+ y((p_2-3\mu)x^2-y^2)\right]\\
\dot{y} &=& \tilde{\lambda} y  +\beta \left[y(p_1 x^2+p_3y^2)+ x(x^2+(p_2+3\mu)y^2)\right]
\end{array}  \right.$ \quad
with $\mu <-\frac{1}{3}$,\quad $\beta=\pm 1$;
\smallbreak

  \item\label{Ana2}
    $\left\{\begin{array}{lcl}
\dot{x} &=& \tilde{\lambda} x  + x(p_1x^2+ p_3y^2)+y((p_2-3\alpha\mu)x^2-\alpha y^2)\\
\dot{y} &=& \tilde{\lambda} y  + y(p_1 x^2+p_3y^2)+ x(\alpha x^2+(p_2+3\alpha\mu)y^2)
\end{array}  \right.$\quad
with $\begin{array}{ll}\alpha =\pm 1&\\ \mu >-\frac{1}{3}&\mu \neq\frac{1}{3};\end{array}$
\smallbreak

  \item\label{Ana3}
  $\left\{\begin{array}{lcl}
\dot{x} &=& \tilde{\lambda} x  +\beta  \left[ x(p_1x^2+p_3y^2)+y((p_2-3\mu)x^2+ y^2)\right]\\
\dot{y} &=& \tilde{\lambda} y  +\beta  \left[ y(p_1 x^2+p_3y^2)+ x(x^2+(p_2+3\mu)y^2)\right]
\end{array}  \right.$\quad
with  $\beta=\pm 1$;
\smallbreak

  \item\label{Ana4}
  $\left\{\begin{array}{lcl}
\dot{x} &=& \tilde{\lambda} x  +  x(p_1x^2+p_3y^2)+y((p_2-3\alpha)x^2-\alpha y^2)\\
\dot{y} &=& \tilde{\lambda} y  + y(p_1 x^2+p_3y^2)+x((p_2+3\alpha)y^2)
\end{array}  \right.$\quad
with $\alpha =\pm 1$;\smallbreak

  \item\label{Ana5}
  $\left\{\begin{array}{lcl}
\dot{x} &=& \tilde{\lambda} x + x(p_1x^2+p_3y^2)+y((p_2-3\alpha)x^2+\alpha y^2)\\
\dot{y} &=& \tilde{\lambda} y + y(p_1 x^2y+p_3y^2)+x(p_2+3\alpha)y^2
\end{array}  \right.$\quad
with $\alpha =\pm 1$;\smallbreak

  \item\label{Ana6}
  $\left\{\begin{array}{lcl}
\dot{x} &=& \tilde{\lambda} x  + x(p_1x^2+p_3y^2)+y((p_2-\alpha)x^2-\alpha y^2)\\
\dot{y} &=& \tilde{\lambda} y  + y(p_1 x^2+p_3y^2)+x(\alpha x^2+(p_2+\alpha)y^2)
\end{array}  \right.$\quad
with $\alpha =\pm 1$;\smallbreak

  \item\label{Ana7}
  $\left\{\begin{array}{lcl}
\dot{x} &=& \tilde{\lambda} x  +  x(p_1x^2+p_3y^2)+y(p_2-3\alpha)x^2\\
\dot{y} &=& \tilde{\lambda} y + y(p_1 x^2+p_3y^2) +x(p_2+3\alpha)y^2
\end{array}  \right.$\quad
with $\alpha =\pm 1$;\smallbreak

  \item\label{Ana8}
 $\left\{\begin{array}{lcl}
\dot{x} &=& \tilde{\lambda} x  +\beta \left[ x((p_1-1) x^2+p_3y^2)+p_2x^2y\right]\\
\dot{y} &=& \tilde{\lambda} y  +\beta \left[y((p_1+3) x^2+p_3y^2)+p_2 xy^2\right]
\end{array}  \right.$\quad
with  $\beta=\pm 1$;\smallbreak

  \item\label{Ana9}
  $\left\{\begin{array}{lcl}
\dot{x} &=& \tilde{\lambda} x  + x(p_1x^2+p_3y^2)+p_2 x^2y\\
\dot{y} &=& \tilde{\lambda} y  + y(p_1 x^2+p_3y^2)+x(\alpha x^2+p_2y^2)
\end{array}  \right.$\quad
with $\alpha =\pm 1$;\smallbreak

  \item\label{Ana10}
  $\left\{\begin{array}{lcl}
\dot{x} &=& \tilde{\lambda} x  + x(p_1x^2+p_3y^2)+p_2x^2y\\
\dot{y} &=& \tilde{\lambda} y  + y(p_1 x^2+p_3y^2)+p_2xy^2 .
\end{array}  \right.$
\end{enumerate}
For $\lambda>0$ these canonical forms give rise to a minimum of 30 qualitatively different phase portraits depicted in Figures~\ref{fig: caseI}--\ref{fig: caseX}.
\end{proposition}

Each normal form in Proposition~\ref{prop:deg3} corresponds to a single expression for $\GG \xy=xQ_2\xy-yQ_1\xy$, that does not depend on  $p_1,p_2,p_3$ and may be multiplied by  $\alpha=\pm 1$ or $\beta=\pm 1$.
Hence the expression for $g(\theta)$ is the same in each group, as well as the angular dynamics at infinity, up to a reversal of direction of rotation at infinity.
This information is summarised in Table~\ref{table:deg3}.
However, the expression of $f(\theta)$, governing the radial dynamics near infinity, depends strongly on these parameters, giving rise to qualitatively different global dynamics.

\begin{table}[hhb]
\begin{tabular}{cclc}
normal	&	equilibria	&	angular	&	figure	 \\
form	&	at infinity	&	stability	&		 \\ \hline
\ref{Ana1}	&	8	&	hyperbolic	&	\ref{fig: caseI}	 \\ \hline
\ref{Ana2}	&	0	&	-	&	\ref{fig: casesIiandVI}	 \\ \hline
\ref{Ana3}	&	4	&	hyperbolic	&	\ref{fig: caseIII}	 \\ \hline
\ref{Ana4}	&	2	&	saddle-nodes	&	\ref{fig:casesIVandIX}	 \\ \hline
\ref{Ana5}	&	6	&	4 hyperbolic	&	\ref{fig: caseV}	 \\ 
	&		&	2 saddle-nodes	&		 \\ \hline
\ref{Ana6}	&	0	&	-	&	\ref{fig: casesIiandVI}	 \\ \hline
\ref{Ana7}	&	4  (0, $\pi/2$, $\pi$, $3\pi/2$)	&	 saddle-nodes	&	\ref{fig: caseIII}	 \\ \hline
\ref{Ana8}	&	4  (0, $\pi/2$, $\pi$, $3\pi/2$)	&	2 hyperbolic	&	\ref{fig: caseIII}	 \\
	&		&	2 hyperbolic-like	&		 \\ \hline
\ref{Ana9}	&	2  ($\pi/2$, $3\pi/2$)	&	saddle-nodes	&	\ref{fig:casesIVandIX}	 \\ \hline
\ref{Ana10}	&	$\infty$	&	-	&	\ref{fig: caseX}	 \\ \hline
&&&
\end{tabular}
\caption{Number and angular stability of infinite equilibria for normal forms in Proposition~\ref{prop:deg3}. Hyperbolic-like are weak non-hyperbolic attractors or repellors.
The last column gives the number of the figure that contains all the possible phase portraits for each normal form.
}\label{table:deg3}
\end{table}

Phase portraits for the normal forms in Proposition~\ref{prop:deg3} may  be obtained applying  Theorem~\ref{teo:globalDynamics} to the information on the behaviour at infinity of the homogeneous differential equations studied in \cite{AnnaL1990}.
When \eqref{eqR2contracting} has finitely many infinite equilibria, their angular stability is determined by the nonlinear part, by Proposition~\ref{prop:stability}.
Since we want to preserve the time orientation, when the infinite equilibria are hyperbolic, each phase portrait in \cite{AnnaL1990} gives rise to potentially two different ones 
obtained applying Theorem~\ref{teo:NewGlobalDynamics}, although in some cases they may coincide, as in cases (0), (1) and (3) of Figure~\ref{fig: caseI}.
This figure contains all the  cases with 8 hyperbolic equilibria at infinity
for $\lambda>0$.
Even if some  infinite equilibrium is not hyperbolic, if $Q_1$ and $Q_2$ have no common factor, then 
Lemmas~\ref{lem:factor f}, \ref{lema:parabolicsector} and \ref{lema:hyperbolicsector} imply that the phase portraits are the same.

The procedure outlined above covers all cases with finitely many equilibria at infinity, whose phase portraits are shown in Figures~\ref{fig: caseI}--\ref{fig:casesIVandIX}, when $Q_1$ and $Q_2$ have no common factor, grouped by number of infinite equilibria.
Cases \ref{Ana2} and \ref{Ana6} with no equilibria at infinity are covered by Theorem~\ref{teo:Bende} and shown in Figure~\ref{fig: casesIiandVI}.

In the degenerate case \ref{Ana10} all points in the equator of the Poincar\'e disk are equilibria.
The next lemma describes the dynamics in the finite domain, shown in Figure~\ref{fig: caseX}.

\begin{lemma}\label{lem:caseX}
Let $D=p_1p_3-p_2^2/4$ and $T=p_1+p_3$, then for the normal form \ref{Ana10} and $\tilde\lambda>0$, we have:
\begin{enumerate}
\renewcommand{\theenumi}{(\alph{enumi})}
\renewcommand{\labelenumi}{{\theenumi}}
\item\label{CXD+T+}
if $D>0$ and $T>0$ then there are no finite equilibria and all infinite equilibria are radially attracting;
\item\label{CXD+T-}
if $D>0$ and $T<0$ then there is a closed curve of attracting finite equilibria encircling the origin;
\item\label{CXD-}
if $D<0$ then there are two invariant half-cones containing curves of attracting finite equilibria separated by two half-cones with radially attracting infinite equilibria;
\item\label{CXD0T+}
if $D=0$ and $T>0$  then there are no finite equilibria and all infinite equilibria are radially attracting;
\item\label{CXD0T-}
if $D=0$ and $T<0$ then a diameter contains no attracting finite equilibria and
the two remaining half-planes contain curves of attracting finite equilibria.
\end{enumerate}
\end{lemma}

\begin{proof} 
Equations with normal form \ref{Ana10} satisfy $g(\theta)\equiv 0$, hence all points at infinity are equilibria.
By Propositions ~\ref{lem:equilibria on radii}  and  \ref{prop:stability} the dynamics on the radius associated to $\theta$ is determined by the sign of $f(\theta)$.
A direct computation shows that in this case $f(\theta)=\FF (\cos(\theta),\sin(\theta))$ where 
$\FF \xy=(x^2+y^2)q\xy$ with $q\xy=p_1x^2+p_2xy+p_3y^2$.
The quantities $D$ and $T$ are, respectively, the determinant and the trace of the symmetric matrix that represents the quadratic form $q\xy$.

If both $D$ and $T$ are positive then $q\xy$ is positive definite, hence $f(\theta)>0$ for all $\theta$, establishing \ref{CXD+T+}. Similarly, when $D>0$ and $T<0$ then $f(\theta)<0$ for all $\theta$, hence every radius contains an attracting finite equilibrium, as in \ref{CXD+T-}. In case \ref{CXD-} the matrix representing $q$ has eigenvalues of opposite sign, so $q\xy$ is negative in a cone and each one of its components contains a curve of  attracting finite equilibria.

If $D=0$ then the eigenvalues of  the matrix of $q\xy$ are 0 and   also $T$. 
If $T>0$ then $f(\theta)\ge 0$ for all $\theta$ and \ref{CXD0T+} follows.
When $T<0$ the infinite equilibria on the diameter that is the eigenspace of $0$ are attracting and each radius on the half-planes determined by  this diameter  contains an attracting finite equilibrium, as in \ref{CXD0T-}.
\end{proof}

\newpage

\phantom{lixo}
\newpage

\phantom{lixo}

\paragraph{Acknowledgements:} The authors are grateful to 
R. Prohens and A. Teruel for  fruitful conversations and to F\'abio Scalco Dias for Example~\ref{ex:grau3}.
The first author was partially supported by the Spanish Research project PID2020--113052GB--I00.
The last two authors were partially supported by Centro de Matem\'atica da Universidade do Porto (CMUP), financed by national funds through FCT - Funda\c{c}\~ao para a Ci\^encia e a Tecnologia, under the project UID/00144/2025.

\newpage
\appendix
\section{Calculations for the examples with degree 2}\label{app:deg2}
 Recall that $\GG \xy=xQ_2\xy-yQ_1\xy$ and $\FF \xy=xQ_1\xy+yQ_2\xy$. 
To these correspond $g(\theta)=\GG (\cos\theta,\sin\theta)$ and $f(\theta)=\FF (\cos\theta,\sin\theta)$.
We use Proposition~\ref{lem:equilibria on radii} to both find the infinite equilibria, by solving $g(\theta)=0$, and determine the dynamics, by looking at the sign of $f(\theta_0)$. 

\vspace{.3cm}

\noindent{\bf Form \ref{Date1} }
 It is $\GG\xy=x^3+y^3 =(x+y)(x^2-xy+y^2)=\dfrac{1}{4}(x+y)((x+y)^2+3(x-y)^2)$ and\\ the infinite equilibria appear at $y=-x$, i.e. $\theta_1=-\pi/4$, $\theta_2=3\pi/4$.\\
We have  $\FF\xy=q_1x^3+(q_2+1)x^2y+(q_1-1)xy^2+q_2y^3$, hence $ \FF(x,-x)=2(q_1-q_2-1)x^3$ and $\sign f( \theta_1)=-\sign f(\theta_2) = \sign (q_1-q_2-1)$. 
 There is  at most one finite equilibrium away from the origin.
The following cases appear,  where labelling (A)--(C) corresponds to Figure~\ref{fig:Date1}:\\
(A) $q_1-q_2< 1$ \quad$\Rightarrow$\quad $f(\theta_1)<0$ and $f(\theta_2)>0$  \quad$\Rightarrow$\quad one finite saddle;\\
(B) $q_1-q_2> 1$ \quad$\Rightarrow$\quad $f(\theta_1)>0$ and $f(\theta_2)<0$  \quad$\Rightarrow$\quad one finite attractor;\\
(C) $q_1-q_2= 1$ \quad$\Rightarrow$\quad $f(\theta_1)=f(\theta_2)=0$   \quad$\Rightarrow$\quad no finite equilibria  away from the origin.

\vspace{.3cm}

\noindent{\bf Form \ref{Date2} }
We have
$\GG\xy=x(x^2-3y^2)$   and infinite equilibria occur for $x=\cos \theta= 0$ and $x=\pm \sqrt{3}y$,
i.e. $\theta_k=(2k+1)\pi/6$, $k=0,\ldots,5$.
 Hence, $g'(\theta_0)=-3$, $g'(\theta_1)=+3$, $g'(\theta_2)=-3$ and $g'(\theta_{k+3})=-g'(\theta_k)$.

It is  $\FF\xy=q_1x^3+(q_2+4)x^2y+q_1xy^2+q_2y^3$ so that  $\FF(0,y)=q_2y^3$ and, for $\varepsilon=\pm 1$,  we obtain $\FF(\sqrt{3}\varepsilon y,y)=4(q_2+\sqrt{3}\varepsilon q_1+3) y^3$.
Thus 
 $f(\theta_{k+3})=-f(\theta_k)$.

\begin{figure}[hhh]
  \includegraphics[width=.4\linewidth]{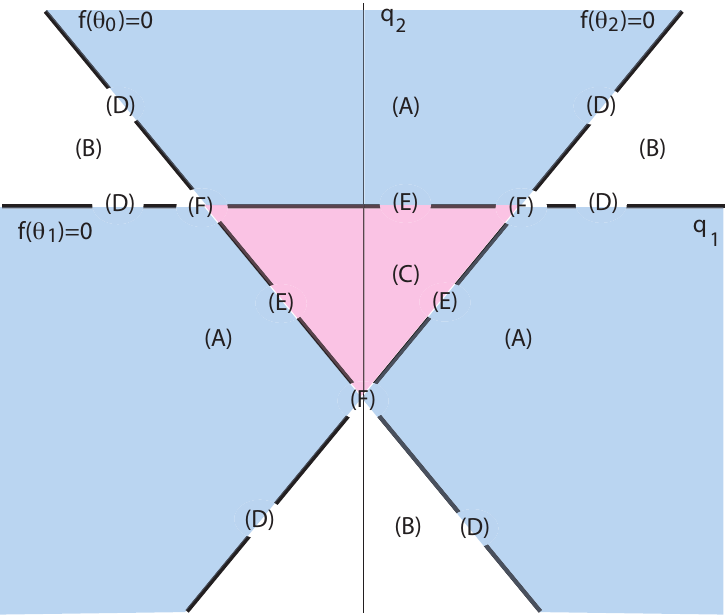}
 \caption{Conditions on the parameters $q_1$ and $q_2$ for the types of phase portrait in normal form \ref{Date2} for $\lambda>0$. 
Regions with the same colour correspond to topologically equivalent phase portraits.
 }\label{fig:Apendice}
 \end{figure}

The following cases occur where the labelling  (A)--(F) refers to the diagrams in Figure~\ref{fig:Date2} for $\lambda>0$ (see Figure~\ref{fig:Apendice}).
 
 Generic cases: three hyperbolic finite equilibria away from the origin.\\
(A) There are two finite saddles, one finite attractor. Conditions:\\
 $q_2>0$ and $q_2>-(\sqrt{3}q_1+3)$ and $q_2>\sqrt{3}q_1-3$ \quad$\Rightarrow$\quad
$f(\theta_0)>0$, $f(\theta_1)>0$, $f(\theta_2)>0$;\\
 $q_2<0$ and $q_2>-(\sqrt{3}q_1+3)$ and $q_2<\sqrt{3}q_1-3$ \quad$\Rightarrow$\quad
$f(\theta_0)>0$, $f(\theta_1)<0$, $f(\theta_2)<0$;\\
 $q_2<0$ and $q_2<-(\sqrt{3}q_1+3)$ and $q_2>\sqrt{3}q_1-3$ \quad$\Rightarrow$\quad
$f(\theta_0)<0$, $f(\theta_1)<0$, $f(\theta_2)>0$.\\
(B) There are  two finite attractors, one finite saddle. Conditions:\\
$q_2>0$ and $q_2>-(\sqrt{3}q_1+3)$ and $q_2<\sqrt{3}q_1-3$ \quad$\Rightarrow$\quad
$f(\theta_0)>0$, $f(\theta_1)>0$, $f(\theta_2)<0$;\\
 $q_2<0$ and $q_2<-(\sqrt{3}q_1+3)$ and $q_2<\sqrt{3}q_1-3$ \quad$\Rightarrow$\quad
$f(\theta_0)<0$, $f(\theta_1)<0$, $f(\theta_2)<0$;\\
$q_2>0$ and $q_2<-(\sqrt{3}q_1+3)$ and $q_2>\sqrt{3}q_1-3$ \quad$\Rightarrow$\quad
$f(\theta_0)<0$, $f(\theta_1)>0$, $f(\theta_2)>0$.\\
(C) There are  three finite saddles. Conditions:\\
 $q_2<0$ and $q_2>-(\sqrt{3}q_1+3)$ and $q_2>\sqrt{3}q_1-3$ \quad$\Rightarrow$\quad
$f(\theta_0)>0$, $f(\theta_1)<0$, $f(\theta_2)>0$.
\smallbreak

Codimension one cases: one pair of non hyperbolic infinite equilibria ($f(\theta_j)=0$ for one $j=0,1,2$) plus two pairs of hyperbolic infinite equilibria. Two hyperbolic finite equilibria away from the origin. Conditions:\\
 $q_2=-(\sqrt{3}q_1+3)$\quad$\Rightarrow$\quad$f(\theta_0)=0$;\\
 $q_2=0$\quad$\Rightarrow$\quad$f(\theta_1)=0$;\\
 $q_2=\sqrt{3}q_1-3$\quad$\Rightarrow$\quad$f(\theta_2)=0$.\\
 (D) two finite equilibria (saddle and attracting node) in consecutive half-cones .
Conditions:\\
$q_2>0$ and $q_2>-(\sqrt{3}q_1+3)$ and $q_2=\sqrt{3}q_1-3$ \quad$\Rightarrow$\quad
$f(\theta_0)>0$, $f(\theta_1)>0$, $f(\theta_2)=0$;\\
$q_2=0$ and $q_2>-(\sqrt{3}q_1+3)$ and $q_2<\sqrt{3}q_1-3$ \quad$\Rightarrow$\quad
$f(\theta_0)>0$, $f(\theta_1)=0$, $f(\theta_2)<0$;\\
$q_2<0$ and $q_2=-(\sqrt{3}q_1+3)$ and $q_2<\sqrt{3}q_1-3$ \quad$\Rightarrow$\quad
$f(\theta_0)=0$, $f(\theta_1)<0$, $f(\theta_2)<0$;\\
$q_2<0$ and $q_2<-(\sqrt{3}q_1+3)$ and $q_2=\sqrt{3}q_1-3$ \quad$\Rightarrow$\quad
$f(\theta_0)<0$, $f(\theta_1)<0$, $f(\theta_2)=0$;\\
$q_2=0$ and $q_2<-(\sqrt{3}q_1+3)$ and $q_2>\sqrt{3}q_1-3$ \quad$\Rightarrow$\quad
$f(\theta_0)<0$, $f(\theta_1)=0$, $f(\theta_2)>0$;\\
$q_2>0$ and $q_2=-(\sqrt{3}q_1+3)$ and $q_2>\sqrt{3}q_1-3$ \quad$\Rightarrow$\quad
$f(\theta_0)=0$, $f(\theta_1)>0$, $f(\theta_2)>0$.\\
(E) two finite equilibria (both saddles) in non consecutive half-cones. Conditions:\\
$q_2=0$ and $q_2>-(\sqrt{3}q_1+3)$ and $q_2>\sqrt{3}q_1-3$ \quad$\Rightarrow$\quad
$f(\theta_0)>0$, $f(\theta_1)=0$, $f(\theta_2)>0$;\\
$q_2<0$ and $q_2>-(\sqrt{3}q_1+3)$ and $q_2=\sqrt{3}q_1-3$ \quad$\Rightarrow$\quad
$f(\theta_0)>0$, $f(\theta_1)<0$, $f(\theta_2)=0$;\\
$q_2<0$ and $q_2=-(\sqrt{3}q_1+3)$ and $q_2>\sqrt{3}q_1-3$ \quad$\Rightarrow$\quad
$f(\theta_0)=0$, $f(\theta_1)<0$, $f(\theta_2)>0$.
\smallbreak

 Codimension two cases: two pairs of non hyperbolic infinite equilibria, one pair of  hyperbolic infinite equilibria.
 Only one hyperbolic finite equilibrium away from the origin.\\
(F) Conditions:\\
$-(\sqrt{3}q_1+3)<0=q_2=\sqrt{3}q_1-3$  \quad$\Rightarrow$\quad
$f(\theta_0)>0$, $f(\theta_1)=0$, $f(\theta_2)=0$;\\
$\pm\sqrt{3}q_1-3=q_2<0$  \quad$\Rightarrow$\quad
$f(\theta_0)=0$, $f(\theta_1)<0$, $f(\theta_2)=0$;\\
$\sqrt{3}q_1-3<0=q_2=-(\sqrt{3}q_1+3)$  \quad$\Rightarrow$\quad
$f(\theta_0)=0$, $f(\theta_1)=0$, $f(\theta_2)>0$.

\vspace{.3cm}

\noindent{\bf Form \ref{Date3} }
When  $\GG \xy=3x^2y$,   infinite equilibria occur at $x=0$ and $y=0$,
i.e. $\theta_k=k\pi/2$, $k=0,\ldots,3$. 
The  infinite equilibria $\theta_1$ and $\theta_3$ are saddle-nodes at infinity,  while $\theta_0$ is a repellor and $\theta_2$ is an attractor. 
There are at most two finite equilibria away from the origin.

We obtain $\FF\xy=q_1x^3+q_2x^2y+(q_1+3)xy^2+q_2y^3$, hence $\FF(x,0)=q_1x^3$ and $\FF(0,y)=q_2y^3$.
Cases, with labelling (A)--(F) as in Figure~\ref{fig:Date3}:\\
(A) $q_1>0$, $q_2\ne 0$  \quad$\Rightarrow$\quad
$f(\theta_0)>0$, $f(\theta_1)\ne 0$  \quad$\Rightarrow$\quad
2 finite equilibria, attractor and saddle-node;\\
(B) $q_1<0$, $q_2\ne 0$  \quad$\Rightarrow$\quad
$f(\theta_0)>0$, $f(\theta_1)=0$  \quad$\Rightarrow$\quad
2 finite equilibria, saddle and saddle-node;\\
(C) $q_1=0$, $q_2\ne 0$  \quad$\Rightarrow$\quad
$f(\theta_0)=0$, $f(\theta_1)\ne 0$  \quad$\Rightarrow$\quad
1 finite equilibrium, saddle-node;\\
(D) $q_1>0$, $q_2=0$  \quad$\Rightarrow$\quad
$f(\theta_0)>0$, $f(\theta_1)=0$  \quad$\Rightarrow$\quad
1 finite equilibrium, attractor;\\
(E) $q_1<0$, $q_2= 0$  \quad$\Rightarrow$\quad
$f(\theta_0)<0$, $f(\theta_1)=0$  \quad$\Rightarrow$\quad
1 finite equilibrium, saddle;\\
(F) $q_1=0$, $q_2=0$  \quad$\Rightarrow$\quad
$f(\theta_0)=0$, $f(\theta_1)=0$  \quad$\Rightarrow$\quad
no finite equilibria away from the origin.

\vspace{.3cm}

\noindent{\bf Form \ref{Date4} }
 When $\GG \xy=x^3$,  equilibria at infinity  occur for $x=0$, i.e.  $\theta_k=k\pi/2$, $k=\pm 1$.
 We have $\FF\xy=q_1x^3+(q_2+1)x^2y+q_1xy^2+q_2y^3$, hence $\FF(0,y)=q_2y^3$  hence $f(\theta_1)$  has the same sign as $q_2$, the opposite sign of  $f(\theta_{-1})$.
The equilibrium at infinity at $\theta_1$ is a non-hyperbolic attractor  at infinity and the one at $\theta_{-1}$ is a non-hyperbolic repellor at infinity.
There is at most one finite equilibrium away from the origin.
Cases, with labelling (A)--(C) as in Figure~\ref{fig:Date4}:\\
(A) $q_2>0$  \quad$\Rightarrow$\quad
$f(\theta_1)>0$  \quad$\Rightarrow$\quad
1 finite equilibrium, saddle;\\
(B) $q_2<0$  \quad$\Rightarrow$\quad
$f(\theta_1)<0$  \quad$\Rightarrow$\quad
1 finite equilibrium, attractor;\\
(C) $q_2=0$  \quad$\Rightarrow$\quad
$f(\theta_1)=0$  \quad$\Rightarrow$\quad
no finite equilibria away from the origin.

\vspace{.3cm}

\noindent{\bf Form \ref{Date5} }
 It is $\GG \xy=0$  so that all points at infinity are equilibria.
 
From  $\FF\xy=q_1x^3+q_2x^2y+q_1xy^2+q_2y^3=(q_1x+q_2y)(x^2+y^2)$, we obtain 
$\sign\FF\xy=\sign(q_1x+q_2y)$.
Finite equilibria apart from the origin satisfy $r=\dfrac{-\lambda}{q_1\cos\theta+q_2\sin\theta}>0$,
or, equivalently, $\lambda+q_1x+q_2y=0$.

\section{Correspondence of the examples with degree 2 with \cite{Artes}}\label{app:correspondence}
The dynamics of vector fields of degree 2 with a star node is described in Theorem 7.2 in the book by Art\'es {\sl et. al} \cite{Artes} for the cases where there are finitely many infinite equilibria.
For the sake of comparison we list here the cases $\alpha_i$, $i=1,\ldots,40$ of that theorem that correspond to each normal form in Proposition~\ref{prop:10vf-2} and to each one of the phase portraits in Subsection~\ref{subsec:deg2}.
From that correspondence we determine whether or not,  at an angle of type $H_0^+$, there is a parabolic sector  delimited by a separatrix that we represented by a dashed line in the phase portraits.

\begin{center}
\begin{tabular}{llclc}case	&	Figure	&	diagram	&	$i$ in $\alpha_i$	&	parabolic	\\	
	&		&		&		&	sector?	\\	\hline
\ref{Date1}	&	\ref{fig:Date1}	&	(A)	&	15--18, 35	&	-	\\	
	&		&	(B)	&	12--14	&	-	\\	
	&		&	(C)	&	29	&	no	\\	\hline
\ref{Date2}	&	\ref{fig:Date2}	&	(A)	&	5--11	&	-	\\	
	&		&	(B)	&	2--4	&	-	\\	
	&		&	(C)	&	1	&	-	\\	
	&		&	(D)	&	27--28	&	no	\\	
	&		&	(E)	&	26	&	no	\\	
	&		&	(F)	&	36	&	no	\\	\hline
\ref{Date3}	&	\ref{fig:Date3}	&	(A)	&	19--20	&	-	\\	
	&		&	(B)	&	21--23	&	-	\\	
	&		&	(C)	&	30	&	no	\\	
	&		&	(D)	&	31--32	&	yes	\\	
	&		&	(E)	&	33--34, 37	&	 no	\\
	&		&	(F)	&	39	&		 no	\\
\hline
\ref{Date4}	&	\ref{fig:Date4}	&	(A)	&	25	&	-	\\	
	&		&	(B)	&	24	&	-	\\	
	&		&	(C)	&	38, 40	&	yes							
\end{tabular}
\end{center}


\noindent
 In case \ref{Date3} (D) both dashed lines exist for some values of the parameter $q_1$.
 \\
 In case \ref{Date4} (C) only one of the dashed lines  may exist, depending on whether it is $\alpha_{38}$ or $\alpha_{40}$.


Thus, the canonical forms in  Proposition~\ref{prop:10vf-2} give rise to a total of 17 qualitatively different phase portraits.

\bigbreak

\end{document}